%% file: sfs.tex
\documentclass[10pt,a4paper,oneside]{article}

\usepackage{latexsym}
\usepackage{amsbsy}
\usepackage{amsmath}
\usepackage{amssymb}
\usepackage{graphics}
\usepackage{epsfig}

\def\QQ{\mathbb{Q}}
\def\RR{\mathbb{R}}
\def\ZZ{\mathbb{Z}}

\newtheorem{thm}{Theorem}[section]
\newtheorem{lem}[thm]{Lemma}
\newtheorem{prop}[thm]{Proposition}

\newenvironment{defin}
      {\vspace{2.5mm}\par\noindent
                \textbf{Definition.}}
             {\vspace{2.5mm}}

\newenvironment{proof}%
      {\par\noindent%
            \textbf{Proof.}}%
           { ~\hfill$\Box$\linebreak}

\newcommand{\fin}{ \!\in\! }
\newcommand{\dv}{ \! : \! }

\title{\large \bfseries TRIANGULATIONS OF SEIFERT FIBRED MANIFOLDS}

\author{\normalsize ALEKSANDAR MIJATOVI\'{C}}

\date{}

\renewenvironment{abstract}
            {\begin{quotation}\noindent\small
                              \textsc{Abstract}.\hspace{0.5mm}}
            {\end{quotation}}

\begin{document}

\maketitle

\begin{abstract}
It is not completely unreasonable to expect that a computable function bounding
the number of Pachner moves needed to change any triangulation of
a given 3-manifold into any other triangulation of the same 3-manifold
exists. In this paper we describe a procedure yielding an explicit
formula for such a function if the 3-manifold in question 
is a Seifert fibred space.
\end{abstract}

\begin{center}
\section{\normalsize \scshape INTRODUCTION}
\label{sec:intro}
\end{center}

\input{intro}

\begin{center}
\end{center}

\begin{center}
\section{\normalsize \scshape SOME NORMAL SURFACE THEORY AND SOME TOPOLOGY}
\label{sec:normal}
\end{center}

\input{somenormal}

\begin{center}
\end{center}

\begin{center}
\section{\normalsize \scshape THE MAIN THEOREM}
\label{sec:main}
\end{center}

\input{main}

\begin{center}
\end{center}

\begin{center}
\section{\normalsize \scshape PACHNER MOVES AND NORMAL SURFACES}
\label{sec:pmns}
\end{center}

\input{pmns}

\begin{center}
\end{center}

\begin{center}
\section{\normalsize \scshape DETECTING SINGULAR FIBRES}
\label{sec:detecting}
\end{center}

\input{fibres}

\begin{center}
\end{center}

\begin{center}
\section{\normalsize \scshape TRIANGULATIONS OF $\mathbf{S^1}$-BUNDLES}
\label{sec:bundle}
\end{center}

\input{bundles1}

\begin{center}
\end{center}

\begin{center}
\section{\normalsize \scshape EXCEPTIONAL MANIFOLDS}
\label{sec:exceptional}
\end{center}

\input{exceptional}

\begin{center}
\end{center}

\begin{center}
\section{\normalsize \scshape CONCLUSION OF THE PROOF OF THEOREM~\ref{thm:main}}
\label{sec:conclusion}
\end{center}

\input{conc}

\nocite{*}
\bibliographystyle{amsplain}
\bibliography{cite}

\noindent \textsc{Department of Pure Mathematics and Mathematical Statistics},
\textsc{Center for Mathematical Sciences},   
\textsc{University of Cambridge}\\
\textsc{Wilberforce Road},
\textsc{Cambridge, CB3 0WB},
\textsc{UK}\\
\textit{E-mail address:} \texttt{a.mijatovic@dpmms.cam.ac.uk}

\end{document}

%% file: intro.tex
\begin{figure}[!hbt]
  \begin{center}
%     \vspace{2cm}  
%     \Huge{(2-3) (1-4)}
    \epsfig{file=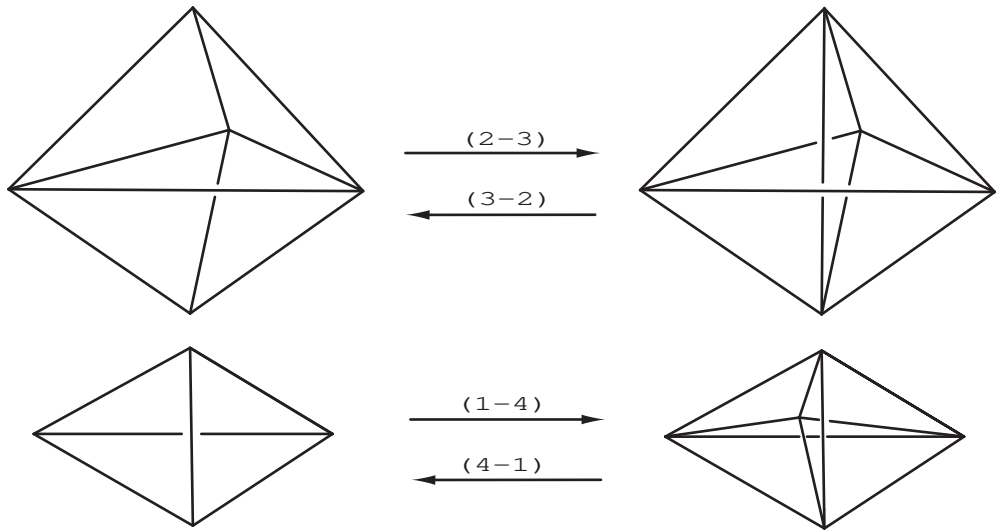}
        \caption{{\small Three dimensional Pachner moves.}}
        \label{fig:3pm}
   \end{center}
\end{figure}

There is a natural way of modifying a triangulation 
$T$
of an
$n$-manifold. Suppose 
$D$
is a combinatorial 
$n$-disc
which is a subcomplex both in this triangulation and in the
boundary of the 
$(n+1)$-simplex 
$\Delta^{n+1}$.
We can change 
$T$
by removing 
$D$
and inserting
$\Delta^{n+1}-\mathrm{int}(D)$.
What we've just described is called a
\textit{Pachner move}. In dimension 
3
there are four possible moves
(see figure~\ref{fig:3pm}).
Note that we can define Pachner moves even if the triangulation
$T$
is non-combinatorial (i.e. simplices 
of
$T$
are not uniquely determined
by their vertices).

Since our aim is to deal with the triangulations of the manifolds 
that are not necessarily 
closed, we need to allow for some additional moves that will
modify the simplicial structure on the boundary
(throughout this paper we will be using the term 
\textit{simplicial structure} as a synonym for a possibly
non-combinatorial triangulation).
The definition of a Pachner move readily generalises to this setting. 
Changing the triangulation of the boundary
by an
$(n-1)$-dimensional Pachner move amounts to gluing onto
(or removing form) our manifold an 
$n$-simplex
$\Delta^{n}$,
which must exist by the definition
of the move. So in dimension 3 we have to use the
three two-dimensional moves (usually referred to as
$(2-2)$,
$(1-3)$
and
$(3-1)$)
that can be implemented by gluing on or shelling a tetrahedron.  

It was proved by Pachner (see~\cite{pachner}) that any two 
triangulations of the same PL 
$n$-manifold are related by a finite sequence of Pachner
moves and simplicial isomorphisms. It is well known 
(see proposition 1.3 in~\cite{mijatov})
that 
a computable function 
%depending only on the 
%number of tetrahedra in the triangulations
%of
%$M$ 
bounding the
length of the sequence from Pachner's theorem in terms 
of the number of tetrahedra
for a fixed 3-manifold 
$M$, 
gives
an algorithm for recognising 
$M$
among all 3-manifolds. The following theorem 
gives an explicit formula for such a bound in case 
$M$
is a Seifert fibred space with a fixed triangulation
on its boundary. \\

\begin{thm}
\label{thm:bullshit}
Let 
$M\rightarrow B$
be a Seifert fibred space with non-empty boundary. Let
$P$
and
$Q$
be two triangulations of 
$M$
that coincide on
$\partial M$
and contain 
$p$
and
$q$
tetrahedra respectively. Then there exists a sequence 
of Pachner moves of length at most
$e^6(10p)+e^6(10q)$
which transforms 
$P$
into a triangulation isomorphic to 
$Q$.
The homeomorphism of 
$M$ 
that realizes the simplicial isomorphism is, when restricted
to
$\partial M$,
equal to the identity on the boundary of
$M$.
\end{thm}

The exponent in the above expression containing the exponential function
$e(x)=2^x$
stands for the composition of the function with itself rather than 
for multiplication. The shameful enormity of the bound can be 
curbed by a more careful choice of subdivisions. The height of
the tower of exponents
can be reduced from 6 to 3, but the complexity of the constructions
involved grows tenfold. Since we are mainly interested in the existence
of an explicit formula, we shall not strive to get the best numbers possible.

The bound in theorem~\ref{thm:bullshit} is clearly computable
and it hence gives a conceptually simple recognition algorithm
for every bounded Seifert fibred space.
It is true that the topology of the Seifert fibred manifolds is not
terribly exciting. They are however ubiquitous
pieces in JSJ-decompositions of 
more interesting 3-manifolds. It is there that our theorem will find
its best application.

%% file: somenormal.tex
The pivotal tool for probing the triangulations of our 3-manifolds
is 
normal surface theory. Many good accounts of it appeared in
the literature (see for example~\cite{hass}, \cite{jaco} or
\cite{bart}). In this section we will state the basic
definitions that 
will allow us to give some known properties
of normal surface. We will then go on to explore
how normal surfaces interact with boundary patterns. The section
will be concluded with a discussion of some well 
known features of incompressible
surfaces in Seifert fibred manifolds. 

We are assuming throughout that all 3-manifolds we are
dealing with are orientable.
Let 
$T$
be a triangulation of a 3-manifold
$M$.
An arc in a 2-simplex of
$T$
is
\textit{normal}
if its ends lie in different 
sides of a 2-simplex. A simple closed
curve in the 2-skeleton of 
$T$
is a \textit{normal curve} if it intersects each 2-simplex of
$T$
in normal arcs. A properly embedded surface 
$F$
in
$M$
is in \textit{normal form} with respect to 
$T$
if it intersects each tetrahedron in 
$T$
in a collection of discs all of whose boundaries are normal
curves consisting of 3 or 4 normal arcs, i.e. triangles and
quadrilaterals. 
A \textit{normal discs} is a triangle or a quadrilateral.
There are precisely seven normal disc types in any tetrahedron
of
$T$.
An isotopy of
$M$
is called a \textit{normal isotopy} with respect to
$T$
if it leaves all simplices of 
$T$
invariant. In particular this implies that it is fixed on the vertices of
$T$.

A normal surface is determined, up to normal isotopy, by the number
of normal disc types in which it meets the tetrahedra of 
$T$.
It therefore defines a vector with
$7t$ 
coordinates. 
Each coordinate represents the number of copies of a normal disc
type that are contained in the surface 
with
$t$
being the number of tetrahedra in
$T$.
It turns out
that there is a certain restricted linear system that such 
a vector is a solution
of. Moreover there is a one to one
correspondence between
the solutions of that restricted linear system and normal 
surface
in
$M$. 
If the sum of two vector solutions of this 
system satisfies the restrictions on the system, then it represents a normal 
surface in
$M$.
On the other hand there is a geometric process called 
\textit{regular alteration} (see figure 2 in~\cite{bart})
which can be carried out
on the normal surfaces representing the summands and which 
yields the normal surface corresponding to the sum. 
It follows directly from the definition of regular alteration
that the Euler characteristic is additive over normal addition.
We can define the 
\textit{weight} 
$w(F)$
of a surface
$F$,
which is transverse to the 1-skeleton of
$T$,
to be the number of points of intersection between 
the surface and the 1-skeleton. Since regular alteration only changes 
the surfaces involved away from the 1-skeleton,
the weight too is additive over normal addition.

A normal surface is called
\textit{fundamental} if the vector corresponding to
it is not a sum of two integral solutions of the linear  
system. 
The solution space of the linear system projects down to a 
compact convex linear cell which is called the
\textit{projective solution space}.
A \textit{vertex surface}
is a connected two-sided normal surface that projects onto a vertex
of the projective solution space (see~\cite{jaco} for a more
detailed description).

The next proposition is proved in~\cite{hass}. It will
be important for us that its proof does not depend on the number 
of equations in the linear system arising from the triangulation
of the manifold.\\

\begin{prop}
\label{prop:hass}
Let 
$M$
be
a compact triangulated 3-manifold containing 
$t$
tetrahedra. Then
each normal coordinate of a vertex surface in
$M$
is bounded above by
$2^{7t}$.
If the normal surface is fundamental, then
$7t2^{7t}$
puts an upper bound on all of its normal coordinates.
\end{prop}

Recall that a properly embedded surface 
$F$
in a 3-manifold 
$M$
is 
\textit{injective} if the homomorphism
$\pi_1(F)\rightarrow \pi_1(M)$,
induced by the inclusion of 
$F$
into 
$M$,
is a monomorphism. A surface 
$F$
is said to be \textit{incompressible} if it satisfies the following conditions:

\begin{list}{$\bullet$}{\itemsep -1mm\topsep-1mm}
\item The surface $F$ does not contain 2-spheres that bound 3-balls
      nor does it
      contain discs
      which are isotopic rel boundary to discs in $\partial M$ and
\item for every disc
      $D$
      in
      $M$
      with 
      $D\cap S=\partial D$
      there is a disc 
      $D'$
      in
      $S$
      with $\partial D=\partial D'$. 
\end{list}
\vspace{1.5mm}

A \textit{horizontal boundary} of an
$I$-bundle over a surface is a part of the boundary
corresponding to the
$\partial I$-bundle. The \textit{vertical boundary} is 
a complement of the horizontal boundary and consists of annuli
that fibre over the bounding circles of the base surface.
It is a well-known fact that a properly embedded one-sided surface
in
$M$
is injective if and only if the horizontal boundary of its regular
neighbourhood is incompressible.
An embedded torus in
$M$
that is incompressible and is not boundary parallel is sometimes
referred to as an \textit{essential torus}.

There is also a relative notion of incompressibility which we
will need to consider. A surface 
$F$
is 
$\partial$-\textit{incompressible} if for each disc
$D$
in
$M$
such that
$\partial D$
splits into two arcs
$\alpha$
and
$\beta$
meeting only at their common endpoints,
with
$D\cap F=\alpha$
and
$D\cap \partial M=\beta$
there is a disc
$D'$
in
$F$
with
$\alpha\subset\partial D'$
and
$\partial D'-\alpha\subset \partial F$.
Such a disc
$D$
is called a 
$\partial$-\textit{compression disc} for 
$F$.
Note that if the manifold
$M$
is irreducible and has incompressible boundary, then we can isotope
$F$
rel
$\alpha$
so that the disc
$D'$
becomes the 
$\partial$-compression disc
$D$.
A properly embedded annulus in 
$M$
that is both incompressible and 
$\partial$-incompressible is called an \textit{essential annulus}. 

Before we can state the main technical results from normal surface theory,
we need to define a very useful concept.

\begin{defin}
A \textit{boundary pattern}
$P$
in a compact 3-manifold
$M$
is a (possibly empty) collection of disjoint
simple closed curves and trivalent graphs embedded in
$\partial M$
such that the surface
$\partial M - P$
is incompressible in
$M$.
\end{defin} 

Boundary patterns usually appear in the context of hierarchies.
In this paper however we will mainly be concerned with a special case when
the pattern
$P$
consists of simple closed curves only. Therefore we will not investigate
the correspondence between 
hierarchies and patterns any further. Notice also that it follows 
from the definition of the pattern that any simple closed curve 
component of
$P$
is homotopically non-trivial in
$\partial M$.

Let
$M$
be a 3-manifold with non-empty
boundary that contains a boundary pattern
$P$.
It follows from the definition that the pattern 
$P$
can be empty 
if and only if
the manifold
$M$
has incompressible boundary. 
Assume now that the pattern 
$P$ 
is not empty.
A subset
of
$M$
is
\textit{pure}
if it has empty intersection with the pattern
$P$.
Most concepts from general 3-manifolds carry over to
3-manifolds with pattern in a very natural way.
For example, a properly embedded surface
$F$
in 3-manifold
$M$
with pattern
$P$
is \textit{P-boundary incompressible}
if for any pure disc
$D$
in
$M$,
such that
$D\cap (\partial M\cup F)=\partial D$
and
$D\cap F$
is a single arc in
$F$,
the arc
$D\cap F$
cuts off a pure disc from
$F$.
Notice that for 
$P=\emptyset$
this notion reduces to 
$\partial$-incompressibility. Also this notion is
well defined only up to an isotopy of
$M$
that is invariant on the pattern
$P$.
In other words we can have two isotopic surfaces in 
$M$,
out of which only one is 
$P$-boundary incompressible.

Let 
$T$
be some triangulation of a 3-manifold 
$M$
with pattern
$P$.
Throughout this paper we will be assuming that the pattern
$P$
lies in the 1-skeleton of
$T$.
This assumption immediately implies that any incompressible
P-boundary incompressible surface in
$M$
can be isotoped into normal form. Also any normal surface
$F$
in
$M$
has a well defined intersection number
$\iota(F)$,
equal to the number of points in
$\partial F\cap P$.
Moreover this  
intersection number is additive over geometric sums of normal 
surfaces.

We say that a surface 
$F$
in a 3-manifold 
$M$
with pattern 
$P$
has \textit{minimal weight} if it can not be isotoped to
a surface with lower weight by
an isotopy that is invariant on the
pattern. In case of
$P=\emptyset$
this reduces to the usual definition of a minimal weight.

An incompressible P-boundary incompressible
surface in 
$M$
of minimal weight has to intersect each triangle in
the 2-skeleton of 
$T$
in normal arcs and possibly some simple closed curves. If
$M$
is irreducible we can remove these simple closed curves by isotopies
in the usual way. The isotoped surface is then in normal form.
The sum
$F=F_1+F_2$
is in \textit{reduced form} if the number of components of
$F_1\cap F_2$
is minimal among all normal surfaces
$F_1'$
and 
$F_2'$
isotopic rel
$P$ 
to 
$F_1$
and 
$F_2$
respectively
such that 
$F=F_1'+F_2'$.
If the pattern
$P$
is empty we get the familiar notion of reduced form which was
used in~\cite{bart}. 
We are now going to define a concept which will be of
some significance to all that follows.

\begin{defin}
A \textit{patch} for the normal sum
$F=F_1+F_2$
is a component of 
$F_1-\mathrm{int}(\mathcal{N}(F_1\cap F_2))$
or
$F_2-\mathrm{int}(\mathcal{N}(F_1\cap F_2))$.
A 
\textit{trivial patch}
is a pure patch which
is topologically a disc and whose boundary intersects 
only one component of the 1-manifold
$F_1\cap F_2$.
\end{defin}

This means that a boundary of a trivial patch is either a single 
simple closed curve in 
$F_1\cap F_2$
or it consists of two arcs: one in
$(\partial M)-P$
and the other in
$F_1\cap F_2$.
The first possibility coincides with what was
called a disc-patch in~\cite{bart}. 
In the absence of pattern the notion of the trivial
patch coincides with what is called a disc-patch
in~\cite{jaco}. The reason for our slightly non-standard
terminology is to avoid the confusion arising from the patches
which are discs but are not disc-patches. There are several
ways in which a disc patch (i.e. a patch which is a disc)
can fail to be trivial, if
$\partial M$
is not empty. Clearly a disc patch which intersects the
pattern
$P$
is non-trivial. A pure disc patch will also be non-trivial
if it intersects
$\partial M$
in more than one arc.

It is not hard to prove that a trivial patch can not have zero 
weight. The argument of lemma 3.3 in~\cite{bart} gives it
to us for all trivial patches bounded by simple closed
curve components of
$F_1\cap F_2$.
In case when our disc patch
is bounded by two arcs,
we can use a simple doubling trick and then apply lemma 3.3
from~\cite{bart} to obtain the desired conclusion.
However there are patches that are topologically discs and have zero
weight. But they must contain more than one component of
$F_1\cap F_2$
in their boundaries.
Now we are finally in the position to state the following lemma.

\begin{lem}
\label{lem:patches}
Let 
$M$
be an irreducible 3-manifold with a (possibly empty) boundary
pattern 
$P$.
Let
$F$
be a minimal weight incompressible P-boundary incompressible
normal surface. If the sum
$F=F_1+F_2$
is in reduced form then each patch is both incompressible and 
P-boundary incompressible and no patch is trivial. Furthermore
if 
$F$
is injective, then each patch has to be injective.
\end{lem}
 
This lemma is a mild generalisation of both lemma 3.6 
in~\cite{bart}
and lemma 6.6 in~\cite{jaco}. Even though patches are not
properly embedded surfaces in
$M$ 
the notions of P-boundary incompressibility and injectivity
can be naturally extended to this setting. Note also that if
$P=\emptyset$,
the manifold 
$M$
has incompressible boundary, the surface 
$F$
is boundary incompressible and so are the patches of
$F=F_1+F_2$.

\begin{proof}
We start by reducing the lemma to the statement that no patch 
of
$F_1+F_2$
is trivial. In case we have a patch which is either compressible
or not injective, we can argue in precisely the same way as in
the proof of lemma 3.6 of~\cite{bart} to obtain a disc patch whose
boundary is a single simple closed curve from 
$F_1\cap F_2$
and is therefore trivial. 

If there is a patch 
$R$
of 
$F_1+F_2$
which has a pure 
boundary compressing disc
$D$,
then we can assume without loss of generality that
$D$
is also a boundary compressing disc for the surface
$F$.
Then the unique arc
$D\cap R=D\cap F$
cuts off a pure disc
$D'$
in
$F$.
If 
$D'$
contains a simple closed curve of
$F_1\cap F_2$,
then we can find a compressible patch of
$F_1+F_2$ 
and we are in the previous case. If there are no simple closed
curves of 
$F_1\cap F_2$
in 
$D'$,
then the edge most arc from
$F_1\cap F_2$
in
$D'$
cuts off a pure disc patch which is clearly
trivial. So it is enough to prove that no 
trivial patch exists.

If there exists a disc patch in 
$F_1+F_2$ 
bounded by a single simple closed curve from
$F_1\cap F_2$,
we can use an identical argument to the one in the proof
of lemma 3.6 in~\cite{bart}
to construct two normal
surfaces
$F'$
and
$S$
such that 
$F=F'+S$.
$F'$
is isotopic to
$F$
and
$S$
is a closed normal surface with
$\chi(S)=w(S)=0$. 
This is clearly a contradiction because no normal surface 
can miss the 1-skeleton. 
Also by our hypothesis the 
surface 
$F$
might be only incompressible
and not necessarily injective. So in order to use the 
argument from~\cite{bart} which shows that every simple
closed curve in
$F_1\cap F_2$,
bounding a trivial patch,
has to be two-sided in both
$F_1$ and 
$F_2$,
we
need to note that 
when 
$M$
equals
$\RR P^3$,
every embedded projective plane in 
$M$
is actually injective.

So now we can assume that 
$F_1+F_2$ 
contains no disc patch disjoint from
$\partial M$.
Let 
$D$
be a trivial patch, lying in
$F_1$
say,
which has least weight among all trivial patches in
$F_1+F_2$.
The intersection
$D\cap F_2$
consists of a unique arc 
$\alpha$
that is a component of
$F_1\cap F_2$.
After regular alteration
$\alpha$
produces two properly embedded arcs in
$F$
one of which
cuts off a pure disc 
$D'$
from 
$F$.
This is because 
$F$
is P-boundary incompressible. Disc
$D'$
is distinct from but might contain the trivial patch
$D$.

We must have
$w(D')= w(D)$.
Otherwise we could isotope 
$F$,
by an isotopy invariant on the pattern,
so that its weight is decreased. 
Since 
$D$
minimises the weight of all trivial patches 
(which is strictly positive)
and
$D'$
must contain at least one such, there is exactly one
trivial patch in 
$D'$
and its weight is equal to 
$w(D)$.
Every other patch in 
$D'$
is topologically a disc whose boundary consists of
4 arcs. Two of them are in 
$\partial M$
and the other two are components of
$F_1\cap F_2$.

If 
$D'$
was itself a patch then we could isotope 
$F_1$
and
$F_2$,
using the definition of the pattern 
$P$
and the irreducibility of 
$M$,
to obtain normal surfaces
$F_1'$
and
$F_2'$
still summing up to 
$F$
but having fewer components of intersection. This contradicts
our assumption on the reduced form of 
$F_1+F_2$.

So 
$D'$
is not itself a trivial patch, but it has to contain one.
Now we can imitate the argument in the proof of lemma 3.6
from~\cite{bart}
to obtain a normal sum
$F=A+F'$
where 
$F'$
is a normal surface isotopic to 
$F$
and 
$A$
is a pure normal annulus of zero weight. Since no normal
surface can live in the complement of the 1-skeleton,
this gives a final contradiction.
\end{proof}

Another crucial fact from normal surface theory, tying up
normal addition with the topological properties of surfaces involved,
is contained in the next theorem. It appeared several times in the literature
in slightly different forms. The version that is of interest to us 
is the following.

\begin{thm}
\label{thm:sum} 
Let 
$M$
be an irreducible 
3-manifold with a possibly empty boundary pattern
$P$. 
Let
$F$
be a least weight normal surface properly embedded in
$M$.
Assume also that 
$F$
is two-sided incompressible 
$P$-boundary incompressible and
$F=F_1+F_2$.
Then 
$F_1$
and
$F_2$
are incompressible and 
$P$-boundary incompressible.
\end{thm}

The proof of this theorem can be obtained by 
using lemma~\ref{lem:patches} 
and following (verbatim) 
the proof of theorem 6.5
in~\cite{jaco}. 
Before we state a further consequence of theorem~\ref{thm:sum},
we need the following simple fact from topology. 

\begin{lem}
\label{lem:annuli}
Let 
$M$
be an irreducible 3-manifold with incompressible boundary. 
Assume also that 
$M$
is neither homeomorphic to
the product 
$S^1\times S^1\times I$,
nor to an
$I$-bundle over a Klein bottle.
Let
$S$
be a toral boundary component of 
$M$.
If 
$A$
and
$B$
are two properly embedded 
incompressible 
$\partial$-incompressible annuli in
$M$
such that
at least one boundary component of both 
annuli lies in 
$S$,
then these bounding simple closed curves must be isotopic
in
$S$,
i.e. they determine the same slope in
$S$.
\end{lem}

\begin{proof}
We start by isotoping the annuli 
$A$
and
$B$
so that their intersection is minimal.
If 
$A\cap B$
is either empty or 
it consists only of essential 
simple closed curves, then the boundary curves must be parallel in
$S$.
So we can assume that 
$\partial A\cap \partial B$
is non-empty. This implies that the boundary slopes of 
$\partial A$
and
$\partial B$
in 
$S$
are distinct. Therefore the complement
$S-(\partial A\cup \partial B)$ 
is a disjoint union of disc.

Since both 
$A$
and 
$B$
are incompressible and 
$\partial$-incompressible, 
the intersection
contains neither contractible simple closed curves nor boundary parallel
arcs in either of the two annuli.
In other words
$A\cap B$
consists of spanning arcs in both annuli. So an
$I$-bundle structure extends from
$A\cup B$
to the the regular neighbourhood
$\mathcal{N}(A\cup B)$.

If the bounding circles of
$A$
lie in distinct components of
$\partial M$,
then the manifold 
$M$
has to be homeomorphic to
$S^1\times S^1\times I$.
This is because each annulus 
$V$
in the vertical boundary of the 
$I$-bundle
$\mathcal{N}(A\cup B)$
cuts off from 
$M$
a 3-ball of the form
$D\times I$
where  
$D$
is one of the discs in
$S-(\partial A\cup \partial B)$.
This 3-ball can not contain 
$A\cup B$
since both annuli are incompressible. We can therefore
extend the product structure over this 3-ball,
thus obtaining an 
$I$-bundle over the torus 
$S$.

If both components of
$\partial A$
live in 
$S$,
then there are two possibilities for the compressible
annulus
$V$.
The dichotomy comes from the discs
$D_1$
and
$D_2$
in the surface 
$S$,
bounded  by the circles of
$\partial V$.
They can either be nested, say 
$D_1\subset \mathrm{int}(D_2)$,
or disjoint.
It is clear that the annulus 
$A$
is disjoint from
$\partial V=\partial D_1\cup \partial D_2$.
The horizontal boundary of the regular neighbourhood
$\mathcal{N}(A\cup B)$
contains an embedded arc, running from 
$\partial D_1$
to 
$\partial A$,
which is disjoint from 
$\partial D_2$.
If the first of the two cases were true, this would
imply that at least one of the boundary components of
$A$
is contained in
$D_2$,
which is clearly a contradiction.
So we must have an embedded 2-sphere
$D_1\cup D_2\cup V$
which bounds a 3-ball 
$D_1\times I$,
disjoint from
$A\cup B$,
like before.
Adjoining all these solid cylinders to 
$\mathcal{N}(A\cup B)$
makes our manifold 
$M$
into an 
$I$-bundle with a single toral boundary component. But this has
to be an 
$I$-bundle
over a closed non-orientable surface of Euler characteristic
zero, i.e. a Klein bottle. This concludes the
proof. 
\end{proof}

A very useful consequence of theorem~\ref{thm:sum} that deals with
orientable surfaces in
$M$
with zero Euler characteristic is contained in proposition~\ref{prop:jaco}.
Most of it is a direct consequence of corollary 6.8 
in~\cite{jaco}.

\begin{prop}
\label{prop:jaco}
Let 
$M$
be an irreducible 3-manifold with incompressible boundary
that contains either an essential torus or an essential
annulus. Let 
$T$
be the triangulation of
$M$.
Then the following
holds:
\newcounter{a}
\begin{list}{(\alph{a})}{\usecounter{a}\parsep 0mm\topsep 1mm}
\item If
      $A$
      is a least weight normal representative (with respect 
      to
      $T$) in an isotopy class
      of an essential torus in
      $M$,
      then every vertex surface in the face of the projective
      solution space that carries 
      $A$
      is an essential torus.
\item Let
      $A$
      be a least weight essential annulus in
      $M$
      that is normal with respect to
      the triangulation
      $T$.
      Then there exists a vertex surface 
      in 
      $M$
      that is an essential annulus. If the boundary
      of
      $M$
      consists of tori only
      and if
      $M$ 
      is neither
      $S^1\times S^1 \times I$
      nor an
      $I$-bundle over a Klein bottle, then the boundary of
      the vertex annulus is parallel to the boundary
      of 
      $A$
      in
      $\partial M$.  
\end{list}
\end{prop}

\begin{proof}
Since vertex surfaces are two-sided by definition, (a) is
just restating corollary 6.8 in~\cite{jaco}. In
(b) at least one of the vertex surfaces in the face of projective 
solution space has to be an annulus. In fact, by theorem~\ref{thm:sum},
it has to be an essential annulus. 

Assume now that 
$M$
has toral boundary. All vertex annuli in the face of the projective
solution space that carries 
$A$
must have their respective boundaries lying in precisely the 
components of
$\partial M$ 
that contain
$\partial A$.
Since all these vertex annuli are essential we can apply 
lemma~\ref{lem:annuli} 
and finish the proof.
\end{proof}

We will conclude this section by a short discussion of
incompressible surfaces in Seifert fibred manifolds.
A surface in a Seifert fibred space 
$M\rightarrow B$
is called \textit{vertical}
if it can be expressed as a union of regular fibres. So the only
possibilities are a torus, a Klein bottle or an annulus. Moebius
band does not come into the picture because the generator of its
fundamental group can not be a regular fibre in 
$M$.
On the other hand the surface that is transverse to all fibres in
$M$
is called \textit{horizontal}. The following proposition 
says roughly that every essential surface in
$M\rightarrow B$
has to be either vertical or horizontal. 

\begin{prop}
\label{prop:hatch}
Let 
$S$
be an incompressible 
$\partial$-incompressible surface in a 3-manifold
$M$.
Then the following holds:
\newcounter{b}
\begin{list}{(\alph{b})}{\usecounter{b}\parsep 0mm\topsep 1mm}
\item Assume further that
$M\rightarrow B$
is an irreducible Seifert fibred space, possibly
containing no singular fibres, over a (perhaps 
non-orientable) compact surface
$B$.
If
$S$
is two-sided and if 
it contains no singular fibres of
$M$,
then it is isotopic to either a vertical surface or a horizontal surface.

\item Let 
$M\rightarrow B$
be an
$S^1$-bundle over a (perhaps non-orientable) compact bounded surface
$B$.
If
$S$
is a one-sided, i.e. non-orientable, 
surface in
$M$,
then it is isotopic to either a vertical
or a horizontal surface.
\end{list}
\end{prop}

It should be noted that (b) of proposition~\ref{prop:hatch} actually
fails if the base surface
$B$
is closed. It is not hard to see that the 
lens space
$L(2n,1)$,
which fibres as an
$S^1$-bundle
over 
$S^2$
with Euler number 
$2n$,
contains an incompressible surface homeomorphic to
a connected sum of 
$n$
projective planes. Such a surface can not be vertical because
of its Euler characteristic, but it can also not be horizontal
since the Euler number of our 
$S^1$-bundle does not vanish. Also notice that the statement 
(b) for two-sided surfaces is already contained in (a).

\begin{proof} If the surface 
$S$
contained a disc component, then, by the definition of incompressibility,
$\partial M$
would have to compress. Since 
$M$
is irreducible, this would make it into a solid torus. The surface
$S$
is then a disjoint collection of compression discs which is
horizontal in any fibration of the solid torus. So from now
on we can assume that 
$S$
contains no disc components.

Part (a) of the proposition is a well-known 
fact about two-sided incompressible 
$\partial$-incompressible surfaces (which
contain no disc components) in Seifert fibred spaces.
Its proof can be found 
in~\cite{jaco1} (theorem VI.34). We shall now prove part
(b).

%Let 
%$C_1,\ldots, C_n$
%be the collection of all singular fibres in 
%$M$.
%If this collection is empty and the base surface 
%$B$
%is closed,
%we take 
%$C_1$
%to be some regular fibre in
%$M$.
%Let
%$M_0$
%be 
%$M$
%with the interiors of small tubular neighbourhoods of
%all 
%$C_i$'s removed.
%We then get an 
%$S^1$-bundle
%$M_0\rightarrow B_0$
%where 
%$B_0$
%is a compact surface obtained from 
%$B$
%by removing 
%$n$
%discs.
%It follows that 
%$\partial B_0$
%can not be empty. Note also that the manifold
%$M_0$
%is still irreducible.

Let 
$M\rightarrow B$
be
an
$S^1$-bundle over a compact bounded surface
$B$.
Choose disjoint arcs in
$B$
whose union decomposes
$B$
into a single disc. Let
$A$
be a union of disjoint vertical annuli
in
$M$
that are the preimages of the collection of arcs under the 
bundle projection.

%We can assume that the surface 
%$S$
%in
%$M$
%is isotoped so that the number of points in the intersection 
%$S\cap (C_1\cup\ldots\cup C_n)$
%is minimal and that it intersects each solid torus 
%$\mathcal{N}(C_1\cup\ldots\cup C_n)$
%in meridional discs.
%Now let
%$S_0$
%be the surface 
%$M_0\cap S$.
%Any compression disc 
%$D$
%for 
%$S_0$
%in
%$M_0$
%is also a compression disc for 
%$S$
%in
%$M$.
%But our hypothesis on
%$S$
%implies that
%$\partial D$
%bounds a disc
%$D'$
%in
%$S$.
%Since 
%$M$
%is irreducible we can isotope 
%$D'$
%rel its boundary
%to 
%$D$.
%If 
%$D$
%was not contained in
%$S_0$,
%this move would decrease the number of points in
%$S\cap (C_1\cup\ldots\cup C_n)$.
%So 
%$S_0$
%must be incompressible in
%$M_0$.
%Similarly a boundary compression disc for 
%$S_0$
%in
%$M_0$
%that intersects
%$\partial M$
%has to cut off a disc in
%$S_0$. 
%Here we are using both the irreducibility of
%$M$
%and the incompressibility of its boundary
%(if
%$\partial M$
%is compressible, then 
%$M$
%is a solid torus and
%$S$
%has to be a disc).
%On the other hand a
%$\partial$-compression disc
%for 
%$S_0$
%in
%$M_0$
%that is disjoint from
%$\partial M$
%and does not cut off a disc in
%$S_0$
%can be used directly to construct an isotopy of
%$S$
%in
%$M$
%that reduces the number of points in
%$S\cap (C_1\cup\ldots\cup C_n)$.
%So
%$S_0$
%is both incompressible and 
%$\partial$-incompressible in 
%$M_0$.

Without loss of generality we can assume that the surface 
$S$
was isotoped in
$M$
so that the number of components of
$S\cap A$
is minimal. 
Then there are no simple closed curves in
$S\cap A$
that are homotopically trivial in either 
$S$
or
$A$,
because such curves can be used to reduce the number of
components in
$S\cap A$
(here we are using the fact that every bounded
$S^1$-bundle is irreducible).
Similarly arcs in
$S\cap A$
that are 
$\partial$-parallel
in either 
$A$
or
$S$
do not occur because both surfaces
are 
$\partial$-incompressible
and 
$\partial M$
is incompressible in 
$M$
(otherwise
$M$
would be a solid torus and 
$S$
would have to be a disc, which is not a one-sided surface).

It follows now that 
$S\cap A$
consists only of vertical circles or horizontal arcs, i.e. spanning arcs
in the annular components of
$A$.
Let 
$M_1$
be the solid torus
$M-\mathrm{int}(\mathcal{N}(A))$
and let
$S_1$
be the surface
$M_1\cap S$.
There can be no trivial simple closed curves in
$(\partial M)\cap M_1$
coming from
$\partial S_1$,
because 
$S$
has no disc components. 
%and it intersects each torus
%from
%$\partial\mathcal{N}(C_1\cup\ldots\cup C_n)$
%in essential simple closed curves.
So we can conclude that 
$\partial S_1$
consists either of horizontal or vertical circles in the torus
$\partial M_1$
(a horizontal simple closed curve is the one that intersects 
each fibre in
$\partial M_1$
transversely).

The surface
$S_1$
has to be incompressible in
$M_1$.
Every compression disc
$D$
for 
$S_1$
in
$M_1$
yields a disc
$D'$
in
$S$.
Using irreducibility of
$M$
we could isotope 
$S$
(rel 
$\partial D$)
so that 
$D'$
becomes
$D$.
If
$D'$
were not contained in
$S_1$,
this move would reduce the number of pieces in
$S\cap A$.

\noindent \textbf{Claim.} If
$\partial S_1$
consists of horizontal simple closed curves, then
$S_1$
is a disjoint union of meridional discs in
$M_1$.

It is enough to show that under the hypothesis of the claim
the surface 
$S_1$
has to be 
$\partial$-incompressible in
$M_1$. 
This is because the only connected incompressible 
$\partial$-incompressible surface in a solid torus
is its meridian disc.

Assume to the contrary that 
$S_1$
is 
$\partial$-compressible.
Let 
$D$
be a 
$\partial$-compression disc for 
$S_1$
in
$M_1$.
We will modify 
$D$,
in a thin collar of 
$\partial M_1$,
so that the arc
$D\cap\partial M_1$
lies in an annulus from 
$\partial M\cap \partial M_1$.
This isotoped 
disc therefore lives in 
$M$
and is a 
$\partial$-compression disc for the surface
$S$.
Like before this leads to a contradiction, because 
we can use the isotoped disc to construct an isotopy
in
$M$
that will reduce the number of
components in
$S\cap A$.

So to prove the claim we need to isotope the 
disc
$D$.
Let
$\alpha$
be the embedded arc
$D\cap\partial M_1$,
running between two points from
$\partial S_1$.
Notice that 
$\partial M_1$
is an alternating union of annuli coming from two families:
one is
$M_1\cap \partial M$
and the other one is
$M_1\cap\mathcal{N}(A)$.
We will now isotope 
$\alpha$
into the interior of
$M_1\cap \partial M$.
Assume first that
$\alpha$
is contained in the interior of an annulus from
$M_1\cap\mathcal{N}(A)$.
Then, since 
$\partial S_1$
is horizontal in
$\partial M_1$,
the segments of
$\partial S_1$
in this annulus
have to be spanning arcs. 
$\alpha$
can either run between two distinct spanning arcs,
or it can run around the annulus to hit the single
spanning arc from two different sides. 
In either situation we can
push 
$\alpha$
into the interior of an adjacent annulus from
$M_1\cap \partial M$.
This isotopy can clearly be extended to the 
collar of
$\partial M_1$, 
thus producing the desired
$\partial$-compression disc.
If, on the other hand,
$\alpha$
is not contained in the interior of an annulus
from
$\partial M_1$,
then we must somewhere have the situation as described
by figure~\ref{fig:hv}.

\begin{figure}[!hbt]
  \begin{center}
    \epsfig{file=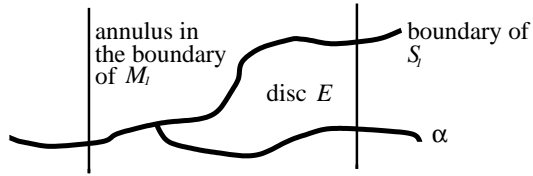}
        \caption{{\small An annulus in $\partial M_1$ containing
                         the disc $E$ with the following property:
                         the intersection
                         $E\cap \alpha$ is a subarc of
                         $\alpha$
                         containing precisely one point from
                         $\partial \alpha$.}}
        \label{fig:hv}
   \end{center}
\end{figure}

We can construct an isotopy  
of 
$\alpha$,
and hence of
$D$,
using the disc 
$E$,
that will reduce the number of points in the 
intersection
$\alpha\cap \partial(M_1\cap \partial M)$.
By repeating this move we arrive at a contradiction
and the claim follows.

It follows directly from the claim that 
$S_1$
with the horizontal boundary can be isotoped 
(rel 
$\partial S_1$)
so that it is horizontal in the fibration of
$M_1$.
This isotopy induces an isotopy of
$S$
in 
$M$
that makes it horizontal.

If 
$\partial S_1$
is vertical,
$S_1$
can only consist of annuli bounded by fibres of
$M_1$.
We can therefore isotope 
$S_1$ (rel 
$\partial S_1$)
so that it becomes vertical. This concludes the proof.
\end{proof}

%% file: main.tex
Let's start by defining precisely the class of Seifert 
fibred 3-manifolds that we shall consider. The manifold
$M$
has to be compact and orientable. It has to fibre into circles
over a possibly non-orientable compact surface 
$B$.
The case when 
$M$
contains no singular fibres will also be considered. 

Since we are only interested in
general Haken 3-manifolds, i.e. the ones that are irreducible and that 
contain an injective surface different from a 2-sphere, 
we need to exclude some Seifert 
fibred spaces. The ones to go first are
$S^1\times S^2$
and the connected sum  
$\RR P^3\#\RR P^3$
since they are not irreducible (the latter manifold fibres as an
$S^1$-bundle over 
a projective plane). 
Another obvious family that we have to disqualify are the lens spaces
including 
the 3-sphere, because they contain no injective surfaces. Those are 
all Seifert fibred spaces over a 2-sphere
with at most 2 singular fibres.

Among the spaces that contain no vertical essential tori are also manifolds
that fibre over a
projective plane
with at most one singular fibre. Since the Klein bottle over the 
orientation reversing curve in 
$\RR P^2$
is not injective, they contain no injective vertical surfaces. There 
are no horizontal surfaces either. The existence of such a surface
would imply that the only singular fibre is of index
$\frac{0}{1}$.
In other words we would be dealing with the connected sum
$\RR P^3\#\RR P^3$
that we've excluded already.
So we need to eliminate all Seifert fibred manifolds with a base surface
$\RR P^2$
that have at most one exceptional fibre because they are not Haken.

The last exceptional class we need to consider are the small Seifert 
fibred spaces. They fibre over 
$S^2$
with three exceptional fibres. It is known that they contain no
separating incompressible surfaces and are therefore Haken
if and only if
$H_1(M;\ZZ)$
is infinite (see for example~\cite{jaco1}, theorem VI.15). 
If this condition is satisfied then every incompressible
surface in our manifold induces a
surface bundle structure over a circle. There are only three small 
Seifert fibred spaces that contain horizontal tori and are therefore
homeomorphic to torus bundles over a circle.
We shall exclude however all small Seifert fibred spaces from further 
consideration, because they contain no essential vertical tori. 

So the 3-manifolds
we are going to consider are precisely all Seifert fibred spaces,
with a possibly empty collection of singular fibres,
that either contain an essential vertical torus or have non-empty boundary.
We can now state the main theorem of this paper.\\

\begin{thm}
\label{thm:main}
Let 
$M\rightarrow B$
be a compact and orientable Seifert fibred 3-manifold that either
contains an essential vertical torus or has non-empty boundary. The base space
$B$
might be non-orientable and the manifold 
$M$
may be an
$S^1$-bundle
containing no singular fibres. Let 
$Y$
be a possibly empty collection of boundary components of 
$M$.
Let
$P$
and
$Q$
be two triangulations 
of
$M$
that contain
$p$
and 
$q$
tetrahedra respectively. Assume also 
that these two triangulations induce isotopic simplicial structures in
$Y$.
Then there exists a sequence of Pachner moves of length at most
$e^6(10p)+e^6(10q)$
which transforms 
$P$
into a triangulation which is isomorphic to 
$Q$.
The homeomorphism of 
$M$
that realizes this simplicial isomorphism maps fibres onto fibres,
does not permute the components of
$\partial M$
and, when restricted to
$Y$,
it is isotopic to the identity on
$Y$.
\end{thm}

The exponent in the expression containing the exponential function
$e(x)=2^x$
stands for the composition of the function with itself and not 
for multiplication.

There are several Seifert fibred manifolds which satisfy
the hypothesis of the theorem but have
more symmetry than the ``generic'' case
(see section~\ref{sec:exceptional}).
What we mean by symmetry is that
such manifolds
contain both vertical and horizontal essential annuli or both
vertical and horizontal essential tori. The
existence of such surfaces
makes it possible to construct homeomorphisms
that don't preserve the fibres. Most of
these manifolds are
$S^1$-bundles with no exceptional fibres.
In this case
the Euler characteristic of the base
surface has to vanish because every horizontal surface
covers it.
If the base space
$B$
has boundary, then
our manifold
is either homeomorphic to
$S^1\times S^1\times I$
or
to the unique orientable
$S^1$-bundle
over a Moebius band. If, on the other hand, the
base surface has no boundary, the manifold
$M$
has to be an
$S^1$-bundle over a torus or over a Klein bottle.
Even though these are infinite classes of manifolds,
there is a unique manifold in each of them that
contains a horizontal torus. This is because
the Euler number
of the these bundles is an integer which has to vanish if
the corresponding Seifert fibred space is to contain a horizontal
surface
(see proposition 4.2 in~\cite{hatch}).
The two manifolds in question are
$S^1\times S^1\times S^1$
and the
$S^1$-bundle over a Klein bottle that is obtained
by identifying two copies of the unique orientable
$S^1$-bundle over a Moebius band via an identity on the
boundary. The triangulations of these manifolds
require a slightly different
technique and will be dealt with in section~\ref{sec:exceptional}.

In the general case the proof of theorem~\ref{thm:main}
is divided into two parts. The
first part is contained in section~\ref{sec:detecting}. In it we
alter the triangulation of
$M$
so that it interacts well with the singular fibres of the fibration.
Theorem~\ref{thm:T2} and proposition~\ref{prop:prod} subdivide
the original triangulation of
$M$
so that a family of vertical solid tori, containing all singular
fibres, are represented by a simplicial subcomplex of the subdivision.
In the second part of the proof we deal with the triangulation
of the 
$S^1$-bundle that is obtained by removing the nicely triangulated 
solid tori. This is described in section~\ref{sec:bundle}.

%% file: pmns.tex
In this section we are going to subdivide, using Pachner
moves, a triangulation
$T$
of a 3-manifold 
$M$
that contains a normal surface 
$F$
of bounded complexity. In other words we are going to alter
$T$
so that its subdivision contains the surface 
$F$
in its 2-skeleton. This procedure will be described by 
lemma~\ref{lem:norm}. It will be applied several
times throughout the construction of the ``canonical'' triangulation of our
Seifert fibred space.\\  

\begin{lem}
\label{lem:norm}
Let 
$M$
be a 3-manifold with a triangulation
$T$
consisting of
$t$
tetrahedra. Assume further that 
$F$
is a properly embedded normal surface in 
$M$
with respect to 
$T$
which contains 
$n$
normal pieces. Than we can obtain a subdivision 
$T_1$
of
$T$,
using less than 
$200nt$
Pachner moves, with the following properties: 
$T_1$
contains the surface 
$F$
in its 2-skeleton and it consists of not more than
$20(n+t)$
tetrahedra. 
\end{lem}

\begin{proof}
Let's start by describing the subdivision 
$T_1$.
Notice that the complement
$M-(T^2\cup F)$
consists of 3-balls,
where 
$T^2$
denotes the 2-skeleton of
$T$. 
The 2-skeleton of
$T_1$
will contain the polyhedron
$T^2\cup F$. 
Its faces are all discs of length at most 6, i.e. their
boundaries consist of at most 6 arcs. The simplicial 
structure on this polyhedron is obtained by coning 
each face from one of the vertices in its boundary.
Now we can define 
$T_1$
to be a union of cones on the boundaries of these 3-balls.
By inspecting all possibilities we see that the number of
triangles in the boundary of any of the above 3-balls is
bounded by 20. Since there are no more than 
$(n+t)$
3-balls, the triangulation
$T_1$
has the desired properties.

The construction of
$T_1$
using Pachner moves will be reminiscent of a similar process
described in section 5 of~\cite{mijatov}.
The main difference here is that the manifold 
$M$
can have a non-empty boundary. If 
$\partial F$ 
is non-trivial, we will also have to change the simplicial structure
on
$\partial M$.
The whole procedure is divided into six stages:

\begin{enumerate}
\item Make a 3-dimensional $(1-3)$ move on each triangle in
      $\partial M$
      and then make a 3-dimensional $(2-2)$ move for each
      1-simplex of 
      $T$
      in
      $\partial M$.
\item Add a vertex into each tetrahedron and each triangle
      of the triangulation
      $T$, 
      that we need to subdivide,
      and then cone.
\item Subdivide the 1-skeleton of
      $T$
      so that it becomes a subcomplex of
      $T_1$,
      and keep the triangulation in the 3-simplices of
      $T$
      coned.
\item Subdivide the 2-skeleton of
      $T$
      to get a subcomplex of
      $T_1$,
      and keep the triangulation in the 3-simplices of
      $T$
      coned.
\item Chop up tetrahedra of
      $T$
      by the appropriate normal pieces coming from
      $F$,
      and
      triangulate the complementary regions by coning them
      from points in their interiors.
\item Shell all tetrahedra that are not contained in any of
      the 3-simplices of 
      $T$.
\end{enumerate}

In the first step we add some tetrahedra to the existing triangulation
of
$M$.
We stop just short of constructing the whole collar on
$\partial M$,
but we do go far enough so that we can later change the original 
simplicial
structure of
$\partial M$
just by using 
$(2-3)$
and
$(1-4)$
moves. In step six we get rid of all redundant 3-simplices that we have 
created in the beginning.

The first step takes 
$10t$
Pachner moves, since there are at most
$4t$
faces in
$\partial M$
and less than 
$6t$
edges. Adding a vertex into 3-simplices in
$T$ 
takes not more than 
$t$
Pachner moves. Adding one into a triangle takes
two Pachner moves. So the second step requires
not more than
$9t$
moves. 

There are at most 
$4n$
vertices of
$T_1$
that are contained in the interiors of 1-simplices of
$T$.
A single vertex can be created on an edge in
$T$
by at most 
$2t$
Pachner moves (see step 2 in section 5 of~\cite{mijatov}).
So step 3 can be accomplished by
$8tn$
moves.

There are at most 
$4n$
normal arcs in all 2-simplices in
$T$.
Each of the complementary discs in each 2-simplex can be triangulated by
at most 4 triangles. So the 2-skeleton of
$T$
will be subdivided by less than
$16n$
triangles. By lemma 4.2 from~\cite{mijatov} this configuration
can be obtained by 
$16n$
two-dimensional Pachner moves. Suspending this process 
gives an upper bound of
$32n$
Pachner moves used in step 4.

In order to glue in the normal discs of
$F$
into the relevant tetrahedra of
$T$
we will need to use the procedure called the changing of cones,
which was described by
lemma 5.1 from~\cite{mijatov}.
We will therefore need to know how many triangles there are
in the disc we are gluing in and also how many 
of them there are in the discs we are changing.
Since we are always gluing in normal pieces, the first number
is bounded above by 2. The second one is bounded by 10, which can be 
seen by inspection. Lemma 5.1 from~\cite{mijatov} now implies that
step 5 can be accomplished by 
$4(2+10)n=48n$
Pachner moves.

It is clear that the triangulation of the ``collar'' that we need to
get rid of in step 6, is shellable. Since Pachner moves that change
the simplicial structure of the boundary are equivalent to elementary
shellings, we only need to count the number of tetrahedra in the
``collar''. 
We know that there are less than 
$16n$
triangles subdividing the original triangulation of
$\partial M$.
In the 3-simplices, that were created by step 1, above the
2-simplices of
$\partial M$
we get not more than
$16n$
tetrahedra.
Each vertex of the subdivision in an edge of 
$T$
lying in 
$\partial M$
gives rise to 2 tetrahedra in the ``collar''.
Since there are not more than 
$4n$
vertices of
$T_1$
contained in the 1-skeleton of
$T$,
the tetrahedra from the second part of step 1 can contain at most
$8n$
3-simplices that need to be shelled. In other words step 6 requires
$24n$
Pachner moves. 
This completes the proof. 
\end{proof}

If the boundary of the manifold 
$M$
is empty, then lemma~\ref{lem:norm} still holds. In fact the procedure 
described in its proof is shortened because steps 1 and 6 become redundant.
Also, if 
$M$
is a submanifold and is triangulated as a subcomplex of the ambient 
triangulation, then the statement of lemma~\ref{lem:norm} remains 
true. A slight modification of step 1 is required, because in this setting
we can alter the triangulation of the boundary of the submanifold using 
$(1-4)$
and
$(2-3)$
moves only.
This situation will arise several times
in sections~\ref{sec:detecting} and~\ref{sec:bundle}.

%% file: fibres.tex
Let 
$M\rightarrow B$
be a Seifert fibred space over a (possibly non-orientable) compact surface
$B$
and let 
$T$
be its triangulation consisting of
$t$
tetrahedra. Our main goal in this
section is to find a subdivision of
$T$
which will support 
in its 2-skeleton an embedded torus around each singular fibre
and will also have some additional properties that will be made precise
later.
We will achieve this by decomposing
$M$
along normal essential tori. 
At the end of the section we will discuss the case when the
surface
$B$
is closed and the manifold 
$M$
contains no singular fibres. In that situation we will subdivide 
$T$
so that its
2-skeleton contains a single embedded torus around some regular fibre.
We are assuming throughout this section that 
the manifold 
$M$
is not homeomorphic 
to one of the five exceptional cases which are described in 
section~\ref{sec:exceptional}. 

Let's assume that our manifold 
$M$
contains an essential vertical torus. We shall now construct a maximal
collection of such tori so that they are pairwise disjoint and
are not topologically parallel. We will also make sure that these
tori have bounded normal complexity. The whole procedure will be 
reminiscent of the strategy used in~\cite{mijatov} to build 
a maximal collection of pairwise disjoint non-parallel normal 
2-spheres. The only major difference in the recursive construction
of 
the family of 
essential tori is that each time we need to produce a new surface
in our collection,
we have to invoke proposition~\ref{prop:jaco}, rather than using
basic normal surface theory, as we did with the 2-spheres.
The estimates for normal complexity of the surfaces
involved here will, however, be identical to the ones used to bound the 
number of normal pieces in the 2-spheres.

We start by taking a vertical incompressible torus
$A'$
that is not boundary parallel in
$M$
and isotoping it so that it becomes both normal with
respect to 
$T$
and so that it has minimal
weight in its isotopy class.
Then proposition~\ref{prop:jaco} gives us an essential
torus
$A$
in
$M$
which is a vertex surface. By (a) of proposition~\ref{prop:hatch}
this torus has to be 
isotopic to either a vertical or a horizontal surface.
In the former case we take the minimal weight representative in
the isotopy class of 
$A$
to be the first
surface 
$A_1$
in our collection. The latter case can not arise if 
$M$
for example has non-empty boundary, because every horizontal
surface
in
$M$
is a branched covering over the base surface 
$B$. 

But if the manifold
$M$
is closed and it contains both a horizontal torus
and a vertical one, 
then its topology 
is very limited. This hypothesis implies that the following 
formula for the Euler characteristic of the base surface 
must hold:
$\chi(B)=\Sigma_i(1-\frac{1}{q_i})$,
where the sum is over all singular fibres
(if there are no singular fibres in
$M$,
the sum on the left should be replaced with zero).
The integers
$q_i$
are just multiplicities, i.e. denominators of the indices
$\frac{p_i}{q_i}$,
of singular fibres in
$M$.
Apart from certain small Seifert fibred spaces that we have already
excluded, in 
the closed case 
we can also get one of the following 3-manifolds:
$S^1\times S^1\times S^1$,
the unique orientable
$S^1$-bundle over a Klein bottle which contains an embedded section,
a fibration over 
$\RR P^2$
with two singular fibres with indices
$\frac{1}{2}$
and
$-\frac{1}{2}$,
and a Seifert fibred space over 
$S^2$
with four singular fibres, two with index
$\frac{1}{2}$
and the other two with index
$-\frac{1}{2}$.
Notice that the second manifold and the fourth manifold 
from the list are homeomorphic.
We will deal with these exceptional cases in section~\ref{sec:exceptional}.

So we can now assume that our vertex torus
$A_1$
is vertical. Its normal complexity is bounded above by 
proposition~\ref{prop:hass}.
We cut 
$M$
open along 
$A_1$
to obtain a manifold
$M_1$
that inherits a polyhedral structure from
$T$.
In~\cite{mijatov} it is described how to adapt normal surface
equations so that their solutions describe closed normal 
surfaces in this polyhedral structure of
$M_1$.
In fact there can be at most 11 distinct disc types in any tetrahedron
of
$T$.

If the manifold
$M_1$
contains an essential vertical torus, then we repeat the above procedure
to obtain a normal torus
$A_2$,
which is again a minimal weight representative
in an isotopy class of a vertex surface. Now 
$A_2$
is necessarily vertical because
$\partial M_1$
is not empty.
Its normal complexity in the initial triangulation
$T$
can be bounded using proposition~\ref{prop:hass}
again. Just like in~\cite{mijatov} it implies that
the number of copies of any given normal disc type 
in the polyhedral structure of the complement
$M_1$
is bounded above by
$2^{11t}$.

It is clear that these upper bounds remain true
in the polyhedral structure of the manifold
$M_2=M_1-\mathrm{int}(\mathcal{N}(A_2))$,
as well as in all subsequent manifolds obtained by the recursive
construction of vertical tori. The Kneser-Haken finiteness
theorem (see~\cite{hempel}) implies that this recursion has to terminate.
Moreover the number of tori in our collection can not exceed 
$20t$.
This is because both Betti numbers
in the formula 
$8t+\beta_1(M;\ZZ)+\beta_1(M;\ZZ_2)$,
bounding the number of disjoint closed
two-sided incompressible surfaces in
$M$,
are smaller than
$6t$.

\begin{figure}[!hbt]
 \begin{center}
    \epsfig{file=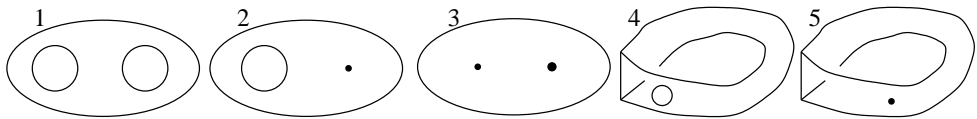}
  \caption{\small The base spaces for all possible bounded Seifert fibred 
                  3-manifolds containing no essential tori.
                  The dots in cases 2 and 3 represent singular fibres
                  and in case 5 it represents  a fibre that may or may not
                  be singular.}  
                  
  \label{fig:base}
 \end{center}
\end{figure}

Now decompose 
$M$
along this disjoint maximal collection of 
non-parallel incompressible tori. 
Since 
$M$
either had an incompressible boundary or it contained an
essential vertical torus,
there are five possible bounded
Seifert fibred spaces that can arise. The respective
base surfaces are depicted in figure~\ref{fig:base}.

Notice also that figure~\ref{fig:base} together with the Kneser-Haken
finiteness theorem implies that any Seifert fibred space, triangulated
by
$t$
tetrahedra,  
can contain at most
$40t$
singular fibres.

Now we form a sub-collection
$\mathcal{A}$
that consists
of some of the essential vertical tori we've constructed. 
We say that an essential torus
belongs to
$\mathcal{A}$
if an only if it bounds on at least one side a Seifert fibred
piece that contains a singular fibre. In other words the
boundaries of cases 2, 3 and possibly 5 from
figure~\ref{fig:base} 
are precisely the tori of
$\mathcal{A}$.
The number of copies of each normal disc type of the triangulation
$T$
in 
$\mathcal{A}$
is bounded above by
$(2\cdot2^{11t})^{20t}$
because at any stage of the construction of
$\mathcal{A}$
there are at most 
$11t$
normal variables (the factor of 2 is a consequence
of the generalisation of normal surface theory to 
the polyhedral setting, see~\cite{mijatov}).\\

We are still 
interested in constructing the family 
$\mathcal{B}$
of compressible tori, one
around each singular fibre of
$M$.
Assume now that we have already triangulated pieces of 
types 2, 3 and 5 from figure~\ref{fig:base} in the
complement of 
$\mathcal{A}$. 
Those are the only complementary pieces of
$\mathcal{A}$ 
that contain singular
fibres. In particular the type 5 piece can not be an
$S^1$-bundle over a Moebius band. 

We proceed by looking for vertical annuli,
that are not 
$\partial$-parallel, in the fibration of the pieces.
If we manage to do that, by finding one such annulus among
vertex surfaces, then it is clear how the horizontal boundary 
of the regular neighbourhood of this
vertex annulus, together with the corresponding
element of
$\mathcal{A}$,
can be used to construct the family 
$\mathcal{B}$
of compressible tori isolating all singular fibres.

Let 
$A$
be a least weight vertical annulus that is not boundary parallel
in a component
$X$
of 
$M-\mathrm{int}(\mathcal{N}(\mathcal{A}))$
that is of
type 2, 3 or 5.
Unless 
$X$
is homeomorphic to an
$I$-bundle
over a Klein bottle, i.e. an
$S^1$-bundle over a Moebius band,
then proposition~\ref{prop:jaco} guarantees the existence
of an essential vertex annulus whose boundary coincides with that
of
$A$.
The vertex annulus is therefore vertical and isotopic to
$A$
(but not necessarily rel
$\partial A$).
In a type 5 piece an element of 
$\mathcal{B}$
consists of the horizontal boundary of the regular
neighbourhood of our vertex annulus, together with
an annulus in the boundary of
$X$.
In the pieces of type 2 and 3 a torus in
$\mathcal{B}$
would consist of a single copy of our vertex annulus together 
with an annulus in
$\partial X$.

The exceptional case, when 
$X$
is equivalent to an
$I$-bundle over a Klein bottle, can actually arise. 
A type 3 piece, with both exceptional fibres of index
$\frac{1}{2}$,
is homeomorphic to an
$S^1$-bundle over a Moebius band. This 
$S^1$-bundle structure is, of course, different from
the original Seifert fibration of 
$X$.
Our vertical annulus 
$A$
becomes a horizontal surface in the
$S^1$-bundle structure. It is in fact equal to the 
horizontal boundary of a regular neighbourhood of some
Moebius band section of this 
$S^1$-bundle.
Up to isotopy however, there are only two possible
Moebius band sections of this 
$S^1$-bundle, one for each isotopy class of orientation
reversing simple closed curves in the vertical Klein 
bottle. Clearly, these two isotopy classes of Moebius band
sections correspond to singular fibres in our original 
Seifert fibration of
$X$. 
If we could make either of these two Moebius bands fundamental,
then we could easily construct compressible tori in
$\mathcal{B}$
around each of the singular fibres. One of them would just be
the boundary of the regular neighbourhood of the 
fundamental Moebius band and the other one would be parallel to the
boundary of its complement in
$X$.
Note also that any properly embedded Moebius band in
$X$
has to be both incompressible and 
$\partial$-incompressible.
The latter follows because 
the surface is one-sided,
$\partial X$
is incompressible and 
$X$
is irreducible. The former is true since the bounding
simple closed curve can not be trivial in
the incompressible boundary of
$X$.
Therefore, by (b) of proposition~\ref{prop:hatch}, any Moebius
band in
$X$
is horizontal in the 
$S^1$-bundle structure on 
$X$.

Take 
$A$
to be the smallest weight essential Moebius band in
$X$.
If we had
$A=F+G$,
we could assume that the sum was in reduced form.
Then lemma~\ref{lem:patches} implies that no
patch is trivial
(in this case boundary pattern 
is empty). 
In other words
$F$
and
$G$
are different from a disc, a 2-sphere and, since
$X$
is irreducible, also from a projective plane.
So we must have
$\chi(F)=\chi(G)=0$.

If
$G$
is closed, then 
$F$
must have the same boundary as
$A$
which consists of a single simple closed curve. That makes
$F$
a Moebius band. 
But this is a contradiction since
the weight of 
$F$
is smaller than that of
$A$.
So both
$F$
and 
$G$
must have boundary.
The same proof tells us that neither of the surfaces
can be a Moebius band. But if they are both annuli, then,
since there are no trivial patches in
$A$,
the
sum must have at least two boundary components which is
again a contradiction.
Therefore 
$A$
has to be fundamental. 
We are now ready to prove the following proposition.

\begin{prop}
\label{prop:vert}
Let
$M\rightarrow B$
be a Seifert fibred space over a possibly non-orientable
compact surface 
$B$
and let 
$T$
be its triangulation consisting of 
$t$
tetrahedra. Assume also that
$M$
contains at least one singular fibre and is either
bounded or it contains an embedded essential 
vertical torus. 
Then there is a subdivision 
$T_1$
of 
$T$
containing the family
$\mathcal{B}$
of compressible tori, one around each singular fibre
of
$M$,
in its 2-skeleton. 
Furthermore
$T_1$
can be obtained from 
$T$
by making less than 
$2^{(2^{400t^2})}$
Pachner moves.
This also gives a bound on 
the number of 3-simplices in
$T_1$.
\end{prop}

\begin{proof}
The construction of
$T_1$
will follow the same lines as the construction of subdivision
$S$
in section 5 of~\cite{mijatov}. The main technical tool
for subdividing 
$T$
with Pachner moves is the procedure described by lemma~\ref{lem:norm}.

Now we can define the subdivision
$T_1$
and then construct it using Pachner moves 
starting from the triangulation
$T$.
Since a single torus in
$\mathcal{A}$
can be in the boundary of two adjacent pieces from
figure~\ref{fig:base},
we need to take doubles of all normal surfaces
contained in
$\mathcal{A}$. 
We'll denote it by
$2\mathcal{A}$. 
The subdivision
$T_1$
will  be obtained in two steps.
In the first step we make sure that the 2-skeleton 
contains the normal surface
$2\mathcal{A}$.
If the family 
$\mathcal{A}$
is empty, then the manifold
$M$
has to be one of the spaces of types 2, 3 or 5 from
figure~\ref{fig:base}. In that case we take
$\mathcal{A}$
to be a single normal torus, parallel to some boundary 
component of
$M$. 
Then we do the same as above.

We know that 
$2\cdot 5t(2\cdot2^{11t})^{20t}$
is an upper bound on the number of normal pieces in
$2\mathcal{A}$. 
Lemma~\ref{lem:norm} implies that this first subdivision 
can be obtained by
$2^{300t^2}$
Pachner moves. This is also a bound on the number of tetrahedra in it.

%build a collar on
%one of the boundary components of
%$M$ and triangulate it so that it contains
%no more than
%$4t$
%tetrahedra. This can be implemented by 
%$3t$
%Pachner moves. We take the boundary component of the collar
%that is contained in the interior of
%$M$
%to be the only element of
%$\mathcal{A}$.
%Notice that the Seifert fibred piece it bounds is triangulate 
%by 
%$t$
%3-simplices.
%
%It is clear that 
%$M-(T^2\cup 2\mathcal{A})$
%consists of 3-balls only. We need to chop these 3-balls
%up by discs coming from our vertical annuli or Moebius bands.
%In the non-parallelity regions of the complement
%$\Delta-2\mathcal{A}$,
%where 
%$\Delta$
%is a 3-simplex in 
%$T$,
%we need to add in PL discs that are described just before
%the statement of proposition~\ref{prop:vert}. Each such disc can
%be triangulated by at most 16 triangles, because
%their length is at most 18, where the \textit{length} of a disc is
%defined to be equal to the number 
%of arcs that its boundary consists of. The number of these discs
%is controlled by the weight of the surface they belong to.
%In the parallelity regions we either need to add horizontal discs
%that are disjoint from
%$\mathcal{A}$
%or vertical discs, but certainly not both. The latter ones arise
%either as parts of boundaries of amalgams in 
%$X$
%or as parts of compression discs when 
%$X$
%is an 
%$in S^1$-bundle over a Moebius band. In both cases the lengths of
%the discs involved are at most 4. So they can be triangulated by
%at most 2 triangle.

In the second step we subdivide further, so that 
singular fibres are isolated in all relevant pieces in
$M-2\mathcal{A}$.
Let
$\mathcal{C}$
be the collection of two parallel 
copies of all vertical annuli constructed 
previously,
in all components of
$M-\mathrm{int}(\mathcal{N}(2\mathcal{A}))$
that are homeomorphic to
type 2, 3 or 5 
Seifert fibred pieces from figure~\ref{fig:base}. 
If a piece of the above complement is a special case
of type 3, then the annulus we have in mind is the horizontal 
boundary of the regular neighbourhood of the fundamental 
Moebius band.
The normal complexity of each of these surfaces is at most four times
the normal complexity of a fundamental surface. Since the union of the relevant
components of
$M-2\mathcal{A}$
is triangulated by less than 
$2^{300t^2}$
tetrahedra, lemma~\ref{lem:norm} and proposition~\ref{prop:hass} imply
that
$$200\cdot2^{300t^2}
(4(5\cdot2^{300t^2})(7\cdot2^{300t^2}2^{7\cdot2^{300t^2}}))<
  2^{2^{400t^2}}$$
bounds the number of Pachner moves required to construct the subdivision 
$T_1$.
The same number also bounds the number of tetrahedra in
$T_1$.
\end{proof}

Our next task is to construct a compression disc in each solid torus
that is bounded by an element of
$\mathcal{B}$.
Notice that the triangulations of these solid tori support,
in their 1-skeleta, a pattern 
$P$
which consists of a single simple closed curve that is
isotopic to a regular fibre (take for example a boundary component of
a surface in
$\mathcal{C}$). 

Let
$D$
be a compression disc that has as few intersections with the 
pattern as possible. In other words
$\iota(D)=q$
where 
$\frac{p}{q}$
is the index of the singular fibre in our solid torus. Assume also
that 
$D$
minimises the weight among all such discs. We will prove that
under these circumstances 
$D$
has to be fundamental.

Assume to the contrary that
$D=F+G$
and that the sum is in reduced form. By lemma~\ref{lem:patches} 
there are no 
trivial patches. 
So neither
$F$
nor
$G$
can be a 2-sphere. They can also not be homeomorphic to
a projective plane because we are living in a
solid torus.
The last possibility
is that
$F$
is a disc.
Again the surface
$G$
can not be closed, because there are no trivial patches.
Since
$\iota(D)=\iota(F)+\iota(G)$,
we must have 
$\iota(D)\geq \iota(F)$.
Since
$w(F)<w(D)$
we can conclude that 
$F$
must be parallel to a disc
$F'$
in the boundary of our solid torus. 
Look at the disc regions in
$F'$
that are bounded by the
arcs coming from 
$\partial G$.
By a nice argument, that is attributed to Haken (see 
claim 4.1.1 in~\cite{bart}), we can conclude that 
one of the disc regions in
$F'$
must, after we do all normal alterations
in
$F+G$,
either be bounded by the whole
$\partial D$
or by a single arc in 
$\partial D$.
The first case can not occur because 
$D$
is a compression disc. In the second case we get an arc
in
$F\cap G$
with one end point in the bad corner of the disc 
region. Using this arc and the disc region we can construct
a pure boundary compression disc for
$D$.
One of the discs, obtained by compressing 
$D$,
will have to be both pure and boundary parallel
(this is because 
$D$
has minimal intersection with the pattern).
This means that the two subdiscs in
$D$,
chopped off by our arc from
$F\cap G$,
are either pure or are disjoint from the pure 
$\partial$-parallel disc we've created by compressing 
$D$.
In the first case we get a pure disc patch in
$D$.
In the second case we can
construct a weight reducing isotopy for 
$D$.
Both of these possibilities lead to contradiction.

Now we know that we can find a fundamental normal disc in every solid 
torus component of the complement of
$\mathcal{B}$.
Assume that we've subdivided the triangulation
$T_1$
from proposition~\ref{prop:vert}, so that these discs are contained
in the 2-skeleton. Their complements in the solid tori
are now 3-balls. We want to change the triangulations in these
3-balls so that they become cones on their respective boundaries.
A slight reformulation of the main result in~\cite{mijatov} 
allows us to do precisely that. 

\begin{thm}
\label{thm:3ball}
Let 
$T$
be a triangulation of a 3-ball with
$t$
tetrahedra. Then it can be changed to a cone on the bounding 2-sphere,
without altering the induced triangulation of the boundary, by less than
$at^22^{at^2}$
Pachner moves, where the constant 
$a$
is bounded above by
$6\cdot 10^6$.
\end{thm}

What we are aiming for at this stage is a subdivision
$T_2$
of the triangulation
$T_1$
that will be ``canonical''
around singular fibres. What we mean by that is that 
a torus 
$\tau$
in
$\mathcal{B}$,
isolating a singular fibre of index
$\frac{p}{q}$,
is triangulated in the following way:
the 1-skeleton in 
$\tau$
contains 
$q$
edges coming from a regular fibre 
$P_\tau$
on 
$\tau$
and
$q$
edges coming from the boundary of the compression 
disc
$D_\tau$. 
This decomposes 
$\tau$
into 
$q$
discs of length four. We triangulate each one of them
by 2 triangles (just pick a diagonal). By doing so,
we haven't introduced any new vertices in the boundary of
$D_\tau$. 
We can therefore triangulate it with 
$q$
triangles
by coning from a point in the disc's interior. The 3-ball region of 
$M-(\tau\cup D_\tau)$
we triangulate by coning as well.
We have thus described the subdivision
$T_2$.
Now we need to construct it using Pachner moves.

The subdivision we have so far is a cone in the 3-ball component
of
$M-(\tau\cup D_\tau)$. 
But the triangulations of 
$\tau$
and
$D_\tau$
are not correct. Let's assume that 
the subdivision
contains two parallel copies of
$\tau$
that have identical simplicial structures and that the product region 
$\tau\times I$
between them is triangulated by cones on the boundary in
every 3-ball of the form
$\Delta \times I$,
where 
$\Delta$
is a triangle from 
$\tau$.
Assume further that the polyhedron
$(\partial D_\tau\cup P_\tau)\times I$
is contained in the 2-skeleton of the subdivision. 
Therefore the arcs of
$\partial D_\tau\cup P_\tau$
are already contained in the 1-skeleton of
$\tau$.
But at the moment they contain many edges of the subdivision.
Later we shall simplify their triangulation to obtain the 
``canonical'' subdivision
$T_2$.
All assumptions we've made here can be implemented while building
the compression discs corresponding to the singular fibres. 
The number of
Pachner moves they require 
can be incorporated in the bounds we have so far.

First we will work on the triangulation of the product
$\tau\times I$.
Let 
$\alpha_\tau$
denote one of the components of
$P_\tau-\partial D_\tau $.
The complement of
$(\partial D_\tau\cup \alpha_\tau)\times I$
in
$\tau\times I$
is a 3-ball which is triangulate by a shellable triangulation.
We can therefore expand the cone structure of 
$\Delta \times I$,
for some triangle 
$\Delta$
from
$\tau$,
to the whole complement (in section 4 of~\cite{mijatov} it is
described how elementary shellings relate to Pachner moves).
Further more we can assume that all three discs in
$((\partial D_\tau\cup \alpha_\tau)-(\partial D_\tau\cap \alpha_\tau))\times I$
are triangulated as cones on their boundaries.
Since all this can be achieved by linearly many Pachner moves,
we can assume that our subdivision had these properties all along.

The PL disc 
$\tau-(\partial D_\tau\cup \alpha_\tau)$
in the subdivision has now cones on both sides. We can therefore change 
its triangulation so that it matches the 
simplicial structure 
induced on it 
by the triangulation
$T_2$,
everywhere but along its boundary. This can again be accomplished
by linearly many Pachner moves
(see lemma 4.2 in~\cite{mijatov}). Now we have to go back to simplifying
the triangulation of the graph
$\partial D_\tau\cup \alpha_\tau$.
This boils down to amalgamating relevant pairs of consecutive edges 
to a single edge. The procedure is the same whether the pair of 
edges we want to get rid of is contained in 
$\partial D_\tau$
or in
$\alpha_\tau$.
This amalgamation is equivalent to crushing one of the edges
in the pair, and thus flattening its star in the subdivision.
This procedure is explicitly described in~\cite{mijatov}.
In our setting it will take
$7=12-6+1$
(cf. lemma 4.1 in~\cite{mijatov})
Pachner moves for one pair of edges. 
Now we can prove the main result of this section.

\begin{thm}
\label{thm:T2}
Let the 3-manifold
$M$
satisfy the same assumptions as in proposition~\ref{prop:vert}.
Then there is a subdivision 
$T_2$
of the triangulation 
$T$
with the properties that are described above. Moreover
$T_2$
can be obtained from
$T$
by making less than
$e^4(500t^2)$
Pachner moves,
where 
$e(x)=2^x$.
The number 
of tetrahedra in
$T_2$
is bounded above by
$e^3(500t^2)$.
\end{thm}

\begin{proof}
We start by applying lemma~\ref{lem:norm} to the fundamental 
compression discs
in the disjoint union of solid tori that are bounded by the
surfaces in
$\mathcal{B}$.
The number of tetrahedra in the solid tori is bounded above
by
$f(t)=e^2(400t^2)$.
So the number of normal pieces in all these compression discs
is bounded above by
$n(t)=40t\cdot5f(t)\cdot7f(t)2^{7f(t)}$.
The factor 
$40t$
bounds the number of singular fibres in
$M$.
So by lemma~\ref{lem:norm} we need to make no more than
$200n(t)f(t)$
Pachner moves to get the compression discs into the 
2-skeleton of the subdivision.

Coning the complementary 3-balls in the solid tori will,
by theorem~\ref{thm:3ball}, take not more than 
$a(20(n(t)+f(t))^2e(a(20(n(t)+f(t))^2)$
Pachner moves, where the constant 
$a$
is bounded above by
$6\cdot10^6$.
Notice that after this step, the number of tetrahedra
in the subdivision is still in the order of magnitude of
$e^3(400t^2)$.
In particular the number of triangles in the compression
discs and toral components of 
$\mathcal{B}$
are of the similar size. Since changing the simplicial structure
of these surfaces, as was described above, can be
done by linearly many (in the number of tetrahedra)
Pachner moves, the number
$e^4(500t^2)$
is surely an upper bound on the total number
of Pachner moves required to build the subdivision
$T_2$.
The inequality 
$20(n(t)+f(t))<e^3(500t^2)$
is also obvious. This completes the proof.
\end{proof}

We shall conclude this section by briefly discussing the case
when 
$M$
is a closed 3-manifold
containing no singular fibres. Since we are dealing with 
Haken 3-manifolds, this implies that the base surface
$B$
has non-zero genus. 
We would like to construct a subdivision
$T_2$
of the triangulation
$T$
of
$M$
so that it contains in its 2-skeleton
a compressible vertical torus. Moreover we want the compression
disc of this 
torus to be a part of the 2-skeleton as well. The 1-skeleton
in the torus will contain both two disjoint parallel regular fibres
and the boundary of this disc. The compression disc itself
will be triangulated by two 2-simplices, while the torus will
consist of four triangles. The complementary region of this two-dimensional
polyhedron that is a 3-ball, will be triangulated as a cone on its boundary.  
This completely describes the subdivision
$T_2$.

We can assume that 
$M$
is neither
$S^1\times S^1\times S^1$
nor
an 
$S^1$-bundle over a Klein bottle which is a double
of the unique orientable
$S^1$-bundle over a Moebius band.
Then proposition~\ref{prop:jaco} gives us a normal vertical torus
that is a vertex surface. Now we can apply lemma~\ref{lem:norm}
to two parallel normal copies of this torus. This
requires not more than
$200\cdot5t2^{7t}t=1000t^22^{7t}$
Pachner moves and produces a subdivision with less than
$20(5t2^{7t}+t)<100t^22^{7t}$
tetrahedra.

Now we look at the complement of the regular neighbourhood
of our vertical torus in
$M$.
Since 
$M$
is a Seifert fibred space and
is not homeomorphic to either of the exceptional manifolds,
at least one component of this complement is not homeomorphic to 
either 
$S^1\times S^1\times I$
or to
an 
$I$-bundle over a Klein bottle.
By proposition~\ref{prop:jaco} we can find
an essential vertical annulus in one of the components of the complement. 
Furthermore this annulus is a vertex surface.
We take two normally parallel copies of it and subdivide further
so that they become part of the 2-skeleton. This can be done
by lemma~\ref{lem:norm}. 
The compressible torus, that we want to have in the 2-skeleton
of
$T_2$, 
will be the boundary of the regular neighbourhood of our 
vertical annulus. Since the triangulation of this
regular neighbourhood, which is a solid torus, is shellable,
we do not need to apply theorem~\ref{thm:3ball}
in order to get the cone. Also, since the
triangulation of the normal vertical annulus contains an embedded spanning
arc in its
1-skeleton, we do not need to look for a compression disc, 
because it is already there (just take the disc above that arc in 
the product structure of the regular neighbourhood).
All the manipulations on the above subdivision, needed to produce
the triangulation
$T_2$,
are very similar to what we did before and also require linearly
many Pachner moves. Thus the following proposition follows.

\begin{prop}
\label{prop:prod}
Let 
$M\rightarrow B$
be an 
$S^1$-bundle over a closed (possibly non-orientable) surface
$B$
that is not homeomorphic to any of the exceptional manifolds
from section~\ref{sec:exceptional}.
Let 
$T$
be its triangulation containing 
$t$
tetrahedra.
Then we can construct a subdivision
$T_2$
of
$T$,
described above,
by making less than
$2^{(2^{400t^2})}$
Pachner moves. 
This also bounds the number of tetrahedra in
$T_2$.
\end{prop}

%% file: bundles1.tex
Let
$B$
be a bounded, possibly non-orientable, compact surface with 
negative Euler characteristic and let
$M\rightarrow B$
be an
$S^1$-bundle over it that is triangulated by
$T$.
Since
$\partial B$
is not empty,
there exists a \textit{section}
of the
$S^1$-bundle
$M$.
That is to say that
$M$
contains a properly embedded surface which is horizontal
and it intersects each fibre precisely once. This is because
$M$
can be obtained from a solid torus
$S^1\times D^2$
by identifying pairs of annuli in the boundary
$S^1\times \partial D^2$.
Each identification can be chosen to be either an identity or
a fixed reflection in the
$S^1$
factor. If
$x\fin S^1$
is invariant under this reflection, then the disc
$\{x\}\times D^2$
yields our section. Notice that we can construct a 
connected horizontal
surface in
$M$,
meeting each fibre in precisely
$n$
points,
for any natural number
$n$,
in a similar way.

Let
$Y$
denote a non-empty collection of boundary components in
$M$.
In this section we are going to construct a triangulation of
$M$
that is going to depend on the triangulation
$T|_Y$
of
$Y$
($T$ is the original triangulation of
$M$).
It is however going to be unique up to homeomorphisms of
$M$
that do not permute components of
$\partial M$
and that are fixed (up to isotopy) on the components
of
$Y$.
For example the simplicial structure in the toral components of
$Y$
naturally arises when our
$S^1$-bundle is a submanifold of a Seifert fibred space
and the tori in
$Y$
bound singular fibres.
We can make sure that
$Y$
is not empty by removing a neighbourhood of a regular
fibre and triangulating the complement as was described
at the end of section~\ref{sec:detecting}.
Roughly speaking we are going to
look for a ``canonical'' section of
$M$
and a family of disjoint vertical annuli that are homologically
non-trivial in
$H_2(M,\partial M;\ZZ_2)$
and whose complement is a solid torus. The union of the ``canonical''
section together with these annuli contains all topological information
about
$M$,
since its complement is just a 3-ball. The triangulation we are looking
for is going to be a cone on the boundary of this 3-ball.

\begin{center}
\end{center}

\begin{center}
\subsection{\normalsize \scshape
                     NORMAL SECTIONS}
\label{subsec:normal}
\end{center}

The ``canonical'' section of 
$M\rightarrow B$
will depend on the simplicial structure of the boundary components in
$Y$.
Before we define it and prove that such a surface is always fundamental
(see theorem~\ref{thm:section}), we need a way of 
relating any two sections
whose boundaries coincide in
$Y$,
by a homeomorphisms of
$M$.
It turns out that twists suffice.
A \textit{twist} along a properly embedded annulus or
torus
in a
3-manifold
$M$
is any homeomorphism of
$M$
which is an identity outside a regular neighbourhood
of the surface. We can now prove the following lemma.\\

\begin{lem}
\label{lem:section}
Let  
$M\rightarrow B$
be an
$S^1$-bundle over a possibly non-orientable bounded surface
$B$
with negative Euler characteristic. Let
$Y$
be a subset of
$\partial M$
as above and let
$S_1$
and
$S_2$
be two sections of the bundle 
$M$,
such that the simple closed curves from
$Y\cap S_1$
and
$Y\cap S_2$
are isotopic in
$Y$.
Then there exists a sequence of twists along vertical annuli in 
$M$ that transforms 
$S_1$
into a section which is isotopic to
$S_2$.
Furthermore this composition of twists restricted
to
$Y$
is isotopic to the identity on
$Y$.
\end{lem}

\begin{proof}
Take a collection of disjoint vertical annuli that decompose
our
$S^1$-bundle into a fibred solid torus. Any section intersects each of these
annuli in a single arc. Now it is clear that any section can be isotoped
so that it coincides with any other section outside a regular
neighbourhood of the vertical annuli. The isotopy can be chosen
so that the surface we are isotoping intersects every fibre of 
the bundle precisely once, throughout the process.
After twisting an appropriate number of times
in a neighbourhood
of a vertical annulus, we can obviously make the two sections
coincide. 

This implies that the sections
$S_1$
and
$S_2$
are related by a sequence of twists along vertical annuli. We can
not guarantee that the annuli we have to twist along are
disjoint from the boundary components in
$Y$.
But each such twist is an orientation preserving
homeomorphism of
$M$
which is 
invariant on the fibres of
$M$.
Also their composition has to preserve the homotopy classes
of the meridians
$S_1\cap Y$
of the tori in
$Y$.
We can therefore conclude that the restriction
of the homeomorphism we get in the end is isotopic to
the identity on
$Y$.
\end{proof}

Let
$T_Y$
denote the triangulation of
$Y$
obtained by restricting the triangulation 
$T$
of the 
$S^1$-bundle
$M$
to the boundary components in
$Y$.
We shall now fix a collection 
$\mathcal{P}$
of simple closed curves in
$Y$
with the following properties:

\begin{list}{$\bullet$}{\itemsep -1mm\topsep-1mm}
\item $\mathcal{P}$ 
      contains precisely one simple closed
      curve on each torus in
      $Y$ 
      and each such curve is a union of normal arcs in
      $T_Y$,
\item $\mathcal{P}$ 
      bounds some normal section 
      $S$
      of the 
      $S^1$-bundle 
      $M$, 
      in other words
      $S\cap Y=\mathcal{P}$, 
      and
\item no collection of simple closed curves satisfying the
      above properties has fewer normal arcs in
      $T_Y$.
\end{list}

A collection 
$\mathcal{P}$
as described here always exists, because our 
$S^1$-bundle 
$M$
always contains a section. Having chosen the collection
$\mathcal{P}$
in this way,
we will now prove (theorem~\ref{thm:section}) 
that we can find a fundamental normal section 
$S$
of
$M$
with the property
$S\cap Y=\mathcal{P}$.
This fundamental surface will therefore interact with
the simplicial structure of
$Y$
in a way which is independent of the triangulation of the
submanifold
$M-Y$.
It is also obvious that lemma~\ref{lem:section}
is precisely what's needed to relate any two such sections.

Before we embark on the proof of theorem~\ref{thm:section},
we need to introduce some notation. For any surface
$F$
that is embedded in the 3-manifold
$M$
and is transverse to the 1-skeleton of 
$T$,
we define the 
``$Y$-weight'',
$w_Y(F)$,
to be the number of points in the intersection of 
$F$
with the 1-skeleton of
$T_Y$.
In other words it follows directly 
from the above definition that the normal section 
$S$
that realizes our chosen collection
$\mathcal{P}$,
minimises the 
``$Y$-weight''
among all normal sections of the bundle 
$M$.
 
\begin{thm}
\label{thm:section}
Let  
$M\rightarrow B$
be an
$S^1$-bundle over a possibly non-orientable bounded surface
$B$
that has negative Euler characteristic.
Let
$Y$
be a non-empty collection of boundary components
of 
$M$
that contains the family 
$\mathcal{P}$
as in the above definition. 
Let
$S$
be a normal section of
$M$,
with respect to the triangulation
$T$,
which 
minimises the weight among all normal sections that 
satisfy the following identity:
$S\cap Y=\mathcal{P}$.
Then
$S$
has to be fundamental.
\end{thm}

\begin{proof}
Assume that 
$S$
can be obtained as a normal sum
$S=F+G$
of two connected surfaces 
$F$
and
$G$.
Without loss of generality
we can take the sum to be in reduced form.
We shall see that this assumption leads into contradiction.

Clearly the homology class
$[S]$
is non-trivial in
$H_2(M,\partial M;\ZZ_2)$,
since it intersects a fibre of
$M$
in a single point. Because we are working over the field
$\ZZ_2$,
we have the following identity 
$[F+G]=[F]+[G]$,
where the first sum is a normal sum and the second is just
addition in the vector space
$H_2(M,\partial M;\ZZ_2)$.
This implies that at least one of the summands,
say
$[F]$,
is different from zero. Now compress and 
$\partial$-compress the surface 
$F$ 
as much as possible. Let
$F'$
be the surface obtained by this process, where we've 
discarded all homologically trivial components.
So 
$F'$
must be homologous to 
$F$
and both
incompressible and
$\partial$-incompressible.
Proposition~\ref{prop:hatch}
implies that 
$F'$
is either vertical or horizontal in the bundle structure of
$M$.

If
$F'$
is vertical, then the homology class 
$[G]$
can not be zero. This is because the boundary
homomorphism
$\delta\dv H_2(M,\partial M;\ZZ_2)\rightarrow H_1(\partial M;\ZZ_2)$
would in this case give
$[\partial S]=\delta([S])=\delta([F]+[G])=\delta([F'])=[\partial F']$.
This is a contradiction because on some boundary torus of
$M$ 
we would have a homology class of the meridian equalling a 
homology class of a fibre. 

If, on the other hand, our surface 
$G$
carries non-trivial homology, then we can effectuate the same process
as before, thus obtaining an incompressible 
$\partial$-incompressible 
surface 
$G'$
which is homologous to 
$G$.
So
$G'$
is either horizontal or vertical. If we had both 
$F'$
and
$G'$
vertical, we would get a contradiction like before, because 
on some boundary torus of 
$M$
we would have a sum of fibres being homologous to a meridian. In other
words, we can assume that one of the two surfaces, say 
$F'$,
is horizontal.

So 
$F'$
has to cover the base surface
$B$
and their respective Euler characteristics are related by
$\chi(F')=k_{F'}\chi(B)$,
for some integer
$k_{F'}$.
Since the inequality
$\chi(F)\leq\chi(F')$
holds, we get the following:
$$\chi(F)\leq\chi(F')=k_{F'}\chi(B)=k_{F'}\chi(S)=k_{F'}(\chi(F)+\chi(G))
  \leq k_{F'}\chi(F).$$
The last inequality follows from the next claim and from the fact that
$M$
contains no embedded projective planes.

\noindent\textbf{Claim.} The surface 
$G$
is neither a disc nor a 2-sphere. 

Since 
$\chi(B)$
is strictly smaller than zero by assumption, the 
same must be true for
$\chi(F')$
and thus for
$\chi(F)$.
So the above inequality is in fact an equality with
$k_{F'}=1$,
$\chi(G)=0$,
and
$\chi(F)=\chi(F')$.
The last equality implies that we neither compressed 
nor
$\partial$-compressed
$F$
while creating 
$F'$.
So 
$F'$
is a subsurface of
$F$.
Since the surface
$F$
is connected,
it must be equal to 
$F'$
and is therefore 
a normal section of
$M$.

On the other hand we have
$w_Y(S)=w_Y(F)+w_Y(G)$
and 
$w_Y(S)\leq w_Y(F)$.
The inequality follows from the definition of the family
$\mathcal{P}$.
This gives us 
$w_Y(G)=0$,
or in other words
$G\cap Y=\emptyset$.
This implies that the normal section 
$F$
also realizes the family
$\mathcal{P}$.
But since 
$G$
is a normal surface, its weight 
$w(G)$
is strictly positive. This yields the inequality 
$w(F)<w(S)$
which contradicts our choice of the normal section
$S$.

So all that is left now is to prove the claim. It is enough to
show that the sum
$S=F+G$
contains no trivial patches
(we take boundary pattern to 
be empty). Even though the sum
$S=F+G$
is in reduced form, we can not directly apply lemma~\ref{lem:patches}
because we couldn't isotope 
$S$
freely to its minimal weight representative
with respect to the triangulation
$T$
(the boundary of
$S$
in
$Y$
was fixed). But away from 
$Y$
the normal surface 
$S$
was chosen so that it minimises the weight. So the same 
argument as in the proof of lemma~\ref{lem:patches}
tells us that any potential trivial patch has to have one arc
from its boundary contained in
$Y$. 

Chose now a trivial patch 
$D$
with 
$w_Y(D)$
minimal among all trivial patches in
$S$.
Like in the proof of lemma~\ref{lem:patches},
using the fact that the surface 
$S$
is boundary incompressible, we get a disc
$D'$
in
$S$
bounded by two arcs: one is in 
$Y$
and the other is ``parallel'' to the arc in the boundary 
of the trivial patch 
$D$.
From the definition of the family
$\mathcal{P}$
it follows that 
$w_Y(D')=w_Y(D)$.
If
$w_Y(D)$
was equal to zero, then we would have both arcs
$D\cap Y$
and
$D'\cap Y$
contained in a single triangle of
$T_Y$.
This would imply that there exist two normal arcs,
one in each summand of the sum
$F+G$,
that intersect in more than one point.
This contradicts lemma 2.1 in~\cite{bart}.
So we can assume that 
$w_Y(D)$
is strictly positive.

If 
$D'$
is a disc patch itself and is hence trivial, then we can isotope
$F$
and
$G$
respectively so that they still sum up to the normal surface
$S$,
but have fewer components of intersection. This 
contradicts the reduced form assumption on the sum
$F+G$
and can therefore not occur.
If 
$D'$
is not a disc patch, then it has to contain
precisely one trivial patch with its 
$Y$-weight
equal to
$w_Y(D)$.
Like in lemma~\ref{lem:patches}, we can construct a 
normal annulus 
$A$,
with its boundary contained in
$Y$,
such that 
$S=A+F'$
and
$w_Y(A)=0$.
But such an annulus does not exist. This proves the claim
and the theorem.
\end{proof}

\begin{center}
\subsection{\normalsize \scshape
                     VERTICAL ANNULI}
\label{subsec:vertical}
\end{center}

So far we having found a ``canonical'' section of
$M$.
Now we need to build a collection of
$(1-\chi(B))$
vertical annuli that are going to decompose our
$S^1$-bundle into a single solid torus. 
Clearly this decomposition is going to be highly non-unique
because the mapping class group of the base surface is full 
of non-trivial elements.
So what we want is to define this collection
in such a way, so that any two possible choices are
related by a homeomorphism of
$M$
which is an identity on the boundary 
$\partial M$.
Also we would like our process
to yield annuli that 
interact with the triangulation
of
$Y$
in a prescribed way.

We start by ordering the toral components of
$\partial M$.
We make sure that in this ordering we enumerate all the
components of
$Y$
before the components of
$\partial M-Y$.
Let 
$Y_1$
in
$Y$
be the first torus in this ordering and let
$A_1,\ldots,A_n$
be a collection of disjoint annuli we are trying to describe
(the subscript
$n$
equals
$1-\chi(B)$).
Each annulus
$A_i$
will contain at least one boundary component in
$Y_1$
and 
every other torus in
$\partial M$
will contain precisely one boundary curve of the surface 
$A_1\cup\ldots\cup A_n$.
We also stipulate that the first 
$2g$
annuli (where
$g$
is the genus of the orientable base surface 
$B$)
have both of their boundary components contained in the torus
$Y_1$.
In case the surface 
$B$
is not orientable, its genus 
$g$
is the maximal number of projective plane
summands it contains when expressed as a connected sum. 
In this case the first 
$g$
annuli from our collection have both of their boundary 
components contained in the torus
$Y_1$.
And finally, again if our base 
$B$
is not orientable, we want to make sure that the complement
of the arc
$p(A_1)$
in the surface 
$B$
is an orientable surface.
The map
$p\dv M\rightarrow B$
denotes our usual bundle projection.
Notice that all of the above conditions on the family of vertical annuli
$A_1,\ldots,A_n$
could be described in terms of their projections onto the base surface.

Now we need to discuss how these vertical annuli interact with 
the triangulation
$T_Y$
of the surface
$Y$.
On each component of
$Y$
we already have a 
normal curve from the collection
$\mathcal{P}$.
The boundaries of the vertical annuli are fibres of 
$M$. 
So the algebraic intersection of a curve 
$\gamma$
in
$\mathcal{P}$
and one such boundary component, both lying in the same torus
of
$Y$,
is 
$+1$
(if the orientations are chosen judiciously). We pick a fibre (of the
bundle
$M$)
$\lambda$
in
each component of
$Y$
with the following property:
$\lambda$
consists of the smallest number of normal arcs in
$T_Y$
in the isotopy class of the fibre (lying in the corresponding component
of
$Y$).
We can assume that simple closed curves
$\gamma$
and
$\lambda$
intersect 
transversely. If they intersect in more than one point,
we will modify 
$\lambda$,
without increasing the number of normal arcs in it,
so that in the end the set
$\gamma\cap \lambda$
contains one point only. This is done in the following way.
The complement of
$\gamma$,
in a component of
$Y$
it lies in,
is an annulus. If 
$\lambda$
intersects this annulus in single spanning arc, then we are done.
If not, there must exist an 
outermost inessential subarc (of
$\lambda$) 
in the annulus, that is parallel
to a subarc of 
$\gamma$
(this is because the intersection number of
$\lambda$
and
$\gamma$
is
$+1$).
It follows from the definition of
$\mathcal{P}$
and the choice of 
$\lambda$
that these two subarcs must contain the same number of normal
arcs.
We can thus modify 
$\lambda$,
without increasing its weight,
so that the number of points in
$\gamma\cap \lambda$
goes down by 2. Repeating this procedure as long as possible 
does the job.

So far we have constructed a normal fibre 
$\lambda$
in each component of 
$Y$
which consists of the smallest number of normal arcs in its isotopy class,
with respect to
$T_Y$.
Furthermore these simple closed curves intersect 
the corresponding elements of
$\mathcal{P}$
in a single point. Now we can define the family 
$\mathcal{F}$
of 
$2n$
fibres, lying in the boundary of
$M$,
which are going to bound our family of vertical annuli
$A_1,\ldots,A_n$.
There will be precisely one fibre on
every torus in
$\partial M-Y_1$.
On each torus of
$Y-Y_1$
we take the simple closed curve 
$\lambda$
we constructed above. Since the triangulation of the boundary
components in
$\partial M-Y$
is not fixed, we can afford to define the simple closed curves 
of 
$\mathcal{F}$
on those tori, only up to isotopy. The distinguished isotopy
class is, of course, that of a fibre.

The situation
on the torus
$Y_1$
is as follows.
If the base surface 
$B$
is orientable (respectively non-orientable), then we take
$2g+1-\chi(B)$
(respectively
$g+1-\chi(B)$)
parallel copies of the simple closed curve 
$\lambda$
defined above. 
It is clear that the family 
$\mathcal{F}$
of normal fibres of 
$M$
that we've just defined,
bounds a collection of vertical annuli 
$A_1,\ldots,A_n$
with all the required
properties. 

Now we need to find a way of constructing
these annuli using normal surface theory. The procedure will
be recursive. So assuming that we have already created 
a subcollection 
$A_1,\ldots,A_k$
(for some 
$k$
smaller than 
$n$)
of vertical normal annuli whose boundaries lie in
$\mathcal{F}$
and which satisfy all other requirements, we look at the manifold
$M_k=M-(A_1\cup\ldots\cup A_k)$
which comes naturally equipped with the polyhedral structure from
the original triangulation
$T$
of
$M$.
Because the annulus 
$A_{k+1}$
we are looking for lives in the manifold
$M_k$
and because its boundary is contained in
$\mathcal{F}$,
i.e. it lies in
$M_k\cap\partial M$,
we can set up the normal surface theory using 
this polyhedral structure on
$M_k$.
How to get the equations and how to bound the complexity
of fundamental surfaces in this normal surface theory is explained
in~\cite{mijatov}. 
From the combinatorial point of view, the procedure
of getting vertical annuli is completely analogous to the process
of building the maximal collection of non-parallel normal 2-spheres
in~\cite{mijatov} and is also very similar to the process of defining
topologically non-parallel vertical tori that was used in
section~\ref{sec:detecting}. All we need now is a way of finding
a fundamental vertical annulus in
$M_k$
that will satisfy all the necessary conditions.\\

\begin{prop}
\label{prop:annuli}
Let
$M\rightarrow B$
be a triangulated
$S^1$-bundle over a possibly non-orientable bounded surface
$B$
with negative Euler characteristic
and let 
$A_1,\ldots,A_k$,
for
$1\leq k<n$,
be normal vertical annuli as described above. Fix the appropriate
simple close curves 
$e$
and
$f$
from 
$\mathcal{F}\cap M_k$,
that are not contained in the union
$A_1\cup\ldots\cup A_k$
and
are supposed to bound the next annulus in our collection.
Let 
$A_{k+1}$
be the normal annulus whose boundary 
consists of the fibres 
$e\cup f$
and which minimises the weight, with respect 
to the polyhedral structure on
$M_k$,
among all normal annuli that carry non-trivial elements of
$H_2(M_k,\partial M_k;\ZZ_2)$
and which are bounded by
$e\cup f$. 
Then
$A_{k+1}$
is fundamental in
$M_k$.
\end{prop}

\begin{proof}
Assume to the contrary that 
$A_{k+1}$
can be expressed as a sum
$A_{k+1}=A+F$,
where 
$A$
and
$F$
are connected normal surfaces. Without loss of generality
we can take it to be in reduced form.  

\noindent\textbf{Claim.} There are no trivial patches in
$A_{k+1}$.\\
Just like in the claim we used to prove theorem~\ref{thm:section},
there can obviously be no trivial patches 
in
$A_{k+1}$
that are disjoint from the boundary tori in
$Y\cap M_k$.
Note also that the boundary curves
$\partial A_{k+1}$
in
$Y\cap M_k$
have, by construction,
the same minimal property the simple closed curves of
$\mathcal{P}$
had in the previous subsection. So we can use an identical argument,
to the one that yielded the claim in the proof of theorem~\ref{thm:section},
to finish the proof of this claim.

This implies that both surfaces
$A$
and
$F$
are different from a disc or a 2-sphere. Since our manifold
$M_k$
contains no projective planes their Euler characteristics are
not positive. Since we have
$0=\chi(A_{k+1})=\chi(A)+\chi(F)$,
both Euler characteristics are in fact equal to zero. 

The claim also implies that the patches on the annulus
$A_{k+1}$
are either discs of length four or annuli. In the former case each
of the discs is bounded by two spanning arcs in
$A_{k+1}$
and by two subarcs in the boundary
$\partial A_{k+1}$.
This means that both surfaces 
$A$
and
$F$
have non-empty boundary. The equation
$[A_{k+1}]=[A]+[F]$
in 
$H_2(M_k,\partial M_k;\ZZ_2)$ 
implies that at least one of them, say
$F$,
is homologically non-trivial.
Since 
$F$
is connected and its
Euler characteristic 
is zero, it has to be both incompressible and boundary incompressible
(here we are using irreducibility of
$M_k$
and incompressibility of its boundary).
By proposition~\ref{prop:hatch} it is therefore either
vertical or horizontal in
the 
$S^1$-bundle structure of
$M_k\rightarrow B_k$. 
But it can not be horizontal if 
$k\geq 2$
because it is disjoint from the annuli 
$A_1\cup\ldots\cup A_{k-1}$.
If 
$k$
equals one,
then it would have to cover the base surface 
$B$
and would thus imply
$\chi(B)=0$.
Since this too is a contradiction, the surface
$F$
must be vertical and therefore an annulus. But on the other hand we
must have at least one boundary component of both  
$F$
and
$A$
contained in
$Y$.
Since 
$F$
is a vertical annulus, its boundary consists of fibres that satisfy
$w_Y(F)<w_Y(A_{k+1})$.
It is true that the boundary of
$F$
can intersect a component of
$\mathcal{P}$
more than once, but we could then improve that 
simple closed curve in
$Y$
(by the same process we used
on
$\lambda$ when defining 
$\mathcal{F}$), 
without increasing the number of normal arcs,
so that the intersection contains a single point, and thus
arrive at a contradiction.

From now on we can assume that all the patches are annuli. This,
together with the fact that 
$A_{k+1}$
is vertical, implies that both surfaces
$A$
and 
$F$
are vertical because they
are just unions of vertical 
patches of
$A_{k+1}$.
If both of them were annuli, then their boundaries 
would have to be disjoint. This can not be, because our annulus
$A_{k+1}$
has only two boundary components. So we can safely assume that
$A$
is an annulus and that 
$F$
is a closed vertical surface. This means that
$\partial A$
equals
$\partial A_{k+1}$.
Our choice of the annulus 
$A_{k+1}$
then gives that 
$[A]$
is trivial in 
$H_2(M_k,\partial M_k;\ZZ_2)$ 
and
that 
$[F]$
carries non-zero homology. In other words
$F$
is not separating in
$M_k$. 

Let
$p\dv M_k\rightarrow B_k$
be the bundle projection, where
$B_k$
is the base surface obtained from 
$B$
by cutting it along the arcs 
$p(A_1)\cup\ldots\cup p(A_k)$.
Then the arc
$p(A)$
separates the base space
$B_k$
into two components
$B'$
and
$B''$.
Also there must exist an arc
$\alpha$
in the family
$p(F-(F\cap A))$,
lying in
say
$B'$,
which is non-trivial in
$H_1(B',\partial B';\ZZ_2)$.
We can now construct a vertical homologically non-trivial
annulus in
$M_k$,
with the boundary equal to
$\partial A_{k+1}$,
whose weight is strictly smaller than that of
$A_{k+1}$.
Let
$\beta$
be a simple closed curve in
$B'$
that intersects
$\alpha$
in a single point (see figure~\ref{fig:dual}).
\begin{figure}[!hbt]
 \begin{center}
    \epsfig{file=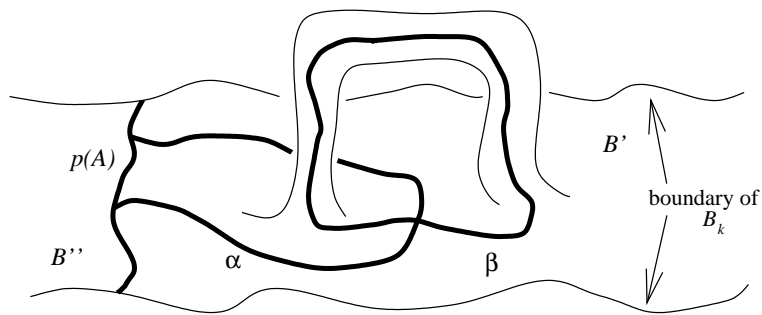}
  \caption{\small The subarc $\alpha$ of the simple closed curve $p(F)$
                  is not homologically trivial in $B'$.}
  \label{fig:dual}
 \end{center}
\end{figure}
Now discarding the subarc of
$p(A)$,
lying between the points of
$\partial \alpha$,
and replacing it with
$\alpha$
yields a properly embedded arc in
$B_k$
which is homologically non-trivial since it intersects
the simple closed curve
$\beta$
in a single point.
The vertical annulus
$A'$
above it represents a non-trivial element in
$H_2(M,\partial M;\ZZ_2)$
and the equation
$\partial A'=\partial A_{k+1}$
holds.
If the weights of the two annuli were the same, then there would exist
a normal isotopy making
$A$
and
$F$
disjoint. Since this can not occur,
$A'$
must have smaller weight than
$A$,
which is again a contradiction.
\end{proof}

The number of vertical annuli in our collection is 
bounded above by the dimension of the vector space
$H_2(M,\partial M;\QQ)$,
which, by Poincar\'{e} duality, equals the first
Betti number
$\beta_1(M;\QQ)$.
This is just a dimension of a vector space
$H_1(M;\QQ)$
that is obtained as a quotient of a space which is 
at most
$6t$
dimensional, where
$t$
is the number of tetrahedra in the 
$S^1$-bundle
$M$.
We therefore have
$n\leq 6t$.

If the base surface of the bundle
$M\rightarrow B$
is orientable, then, using proposition~\ref{prop:annuli},
we can construct our maximal collection 
$A_1,\ldots,A_n$
of vertical normal annuli that has all the required properties.
Furthermore we have control over its normal complexity.
Just like, when bounding the normal complexity of the family
$\mathcal{A}$
of essential vertical tori in section~\ref{sec:detecting},
we can conclude that the number of copies of each normal 
disc type of the triangulation of 
$M$,
that are contained in
$A_1\cup\ldots\cup A_n$,
is bounded above by
$(2\cdot 2^{11t})^{6t}$.

In case the base surface 
$B$
is non-orientable, we need to fulfill one more requirement.
Namely the complement in
$M$
of the first annulus 
$A_1$
in our collection has to be an
$S^1$-bundle
over an orientable surface. 
If we are unlucky, none of the annuli that are recursively 
constructed by proposition~\ref{prop:annuli} will have this property.
If we nevertheless construct the first 
$g$
(where
$g$
denotes the ``non-orientable'' genus of the base surface
$B$)
vertical annuli
$A_1,\ldots,A_g$
in the way described above, using proposition~\ref{prop:annuli},
then their boundaries are contained in the 
torus
$Y_1$.
But more to the point,
the 3-manifold
$M_g=M-(A_1\cup\ldots\cup A_g)$
fibres over a planar, and therefore orientable, connected surface
$B_g$.
Now look at the 
$g$
arcs
$p(A_1),\ldots,p(A_g)$
in the surface
$B$,
where 
$p\dv M\rightarrow B$
is the bundle projection. Their end points are contained in
the circle
$p(Y_1)$.
By concatenating some subcollection of these arcs, we can obviously
construct a single arc 
$\alpha$
such that the complement
$B-\mathrm{int}(\mathcal{N}(\alpha))$   
is a compact orientable surface. 
Also, this construction gives us control over the weight
of the vertical annulus
$p^{-1}(\alpha)$
above the arc 
$\alpha$.
This is because the concatenation of arcs can be implemented
in the 
$S^1$-bundle
$M$
by using some of the annuli in
$Y_1-(A_1\cup\ldots\cup A_g)$
that are pushed slightly into the interior of
$M$.
Such annuli will have zero weight.

So in case when 
$M$
fibres over a non-orientable base space 
$B$,
we choose the first annulus 
$A_1$
to be the smallest weight annulus that satisfies all our
requirements. It follows from what was said above that 
there exists a normal annulus, containing less than 
$(2\cdot 2^{11t})^{6t}$
copies of each normal disc type in the triangulation of
$M$,
that does the job. This is because
$g<n\leq 6t$,
where 
$t$
is the number of tetrahedra in
$M$.
The rest of the vertical annuli in the collection
$A_1,\dots, A_n$
can be obtained as before, using proposition~\ref{prop:annuli}.
The number of copies of a given disc type in the 
surface
$A_1\cup\ldots\cup A_n$,
is certainly bounded 
above by
$(2\cdot 2^{11t})^{12t}$.

Let 
$S$
be a minimal weight normal section of the triangulated
$S^1$-bundle
$M$,
as described by theorem~\ref{thm:section}. Without loss of 
generality we can assume that the intersection 
between the surface
$S$
and our normal vertical annuli, is transversal and
that it misses the 1-skeleton. It is clear that 
the horizontal surface 
$S$
can be isotoped so that it intersects the vertical surface 
$A_1\cup\ldots\cup A_n$
in precisely 
$n$
arcs. But we want this isotopy not to increase the weight of 
either of the two surfaces. 

If there exists a simple closed
curve of intersection, then it has to be inessential in the
annulus it lies in. This is because the generic intersection
contains a spanning arc on the annulus. 
Since
$S$
injects on the
$\pi_1$ level, this simple closed curve 
has to be homotopically trivial in our section 
as well. It therefore bounds a disc in
$S$.
So the innermost 
simple closed curve of intersection
that bounds the lightest disc in one of the surfaces,
has to bound a disc of the same weight in the other surface. 
This is because otherwise we would have a weight reducing isotopy of 
one of the surfaces. 
Moreover,
after we do the exchange, i.e. isotope over the 3-ball 
bounded by the two discs, the surfaces we get are in normal form. 
Repeating this procedure eliminates all simple closed curves of 
intersection between
$S$
and
$A_1\cup\ldots\cup A_n$. 
If there are any inessential arcs of intersection in our vertical
annuli, then their boundaries have to be contained in
$\partial M-Y$.
This is because each fibre in 
$\mathcal {F}$
intersects every simple closed curve in
$\mathcal{P}$
at most once. So using the fact that 
$S$
is 
$\partial$-incompressible, we can isotope it like before, so that 
the intersection 
$S\cap(A_1\cup\ldots\cup A_n)$ 
consists of spanning arcs only. This isotopy will not increase the 
weight of 
$S$. 
The surface we obtain in the end is still in normal form (otherwise we 
could decrease its weight further). 

All we need now is a way of transforming one possible collection
of vertical annuli, that we've constructed, to another. For achieving  
this we will
need the following lemma.

\begin{lem}
\label{lem:twist}
Let
$B$
be a possibly non-orientable bounded surface and let
$\alpha$
and
$\beta$
be two properly embedded arcs in
$B$
with the properties:
the elements of
$H_1(B,\partial B;\ZZ_2)$
carried by arcs
$\alpha$
and
$\beta$
are non-zero,
$\partial \alpha=\partial \beta$,
and
both surfaces
$B_\alpha=B-\mathrm{int}(\mathcal{N}(\alpha))$
and
$B_\beta=B-\mathrm{int}(\mathcal{N}(\beta))$
are orientable.
Then there exists a homeomorphism of
$B$
which is an identity on the boundary
and which maps
$\alpha$
onto
$\beta$.
\end{lem}

\begin{proof}
Let
$B_\alpha'$
and
$B_\beta'$
denote the surfaces
$B_\alpha$
and
$B_\beta$
with discs glued to all of their boundary components.
This means that both
$B_\alpha'$
and
$B_\beta'$
are closed orientable surfaces.
The number of boundary components of
$B_\alpha$
and
$B_\beta$
can differ by at most 1.
But their Euler characteristics have to be equal:
$\chi(B_\alpha)=\chi(B)+1=\chi(B_\beta)$.
On the other hand
the Euler characteristics of
$B_\alpha'$
and
$B_\beta'$
are both even,
so we must have
$\chi(B_\alpha')=\chi(B_\beta')$.
So
$B_\alpha'$
is homeomorphic to
$B_\beta'$.
Moreover
$B_\alpha$
and
$B_\beta$
are connected orientable surfaces with the same number of
boundary components. They are therefore also homeomorphic.

Let
$a$
be one of the two arcs in
$\partial B-(\partial \alpha)$.
Choose an orientation for it.
Notice that
$a$
is contained in
$\partial B_\alpha$
as well as in
$\partial B_\beta$.
This induces an orientation on one boundary
component of
$B_\alpha$
and on a boundary component of
$B_\beta$.
Since both of these surfaces
are orientable, we can extend this
orientation over them. This therefore induces orientations of 
closed surfaces
$B_\alpha'$
and
$B_\beta'$,
and of the simple closed curves
$\partial B_\alpha$
and
$\partial B_\beta$ 
lying in them.
By the classification of closed orientable surfaces
we can find an orientation-preserving homeomorphism
$f'\dv B_\alpha'\rightarrow B_\beta'$
with the following properties:
$f'$
maps homotopically trivial simple closed curves
$\partial B_\alpha$
in
$B_\alpha'$
onto the family
$\partial B_\beta\subset B_\beta'$,
it preserves the orientations of the curves, and acts as
an ``identity'' both on all components of
$\partial B_\alpha$
that are disjoint from
$\alpha$
and on the arcs from
$\partial B_\alpha-(\partial\alpha)$.
The two arcs
$\partial B_\alpha\cap \mathcal{N}(\alpha)$
are mapped onto
$\partial B_\beta\cap \mathcal{N}(\beta)$
by
$f'$,
with their induced orientations preserved.
All this is possible because we can always 
isotope a finite disjoint collection of discs in any surface onto
any other such collection. Once this is established,
everything else follows from the fact that our homeomorphism is
orientation preserving.

It is clear that the homeomorphism
$f'$
induces a homeomorphism
$f\dv B_\alpha\rightarrow B_\beta$
between the two bounded surfaces. It also follows that the map 
$f$
extends to a homeomorphism from
$B=B_\alpha\cup\mathcal{N}(\alpha)$
to the surface
$B=B_\beta\cup\mathcal{N}(\beta)$.
This extension is an identity on the boundary of
$B$.
\end{proof}

Using lemma~\ref{lem:twist} and a well-known fact that any 
homeomorphism of the base surface
$B$,
that is an identity on the boundary of 
$B$,
extends to a homeomorphism of the whole 
$S^1$-bundle
$M\rightarrow B$
which is fixed on
$\partial M$,
tells us that any collection of our vertical annuli
can be obtained from any other collection by such
a homeomorphism.

\begin{center}
\end{center}

\begin{center}
\subsection{\normalsize \scshape
                   THE SIMPLIFIED TRIANGULATION}
\label{subsec:simplified}
\end{center}

Let
$T$
be a triangulation of the 
$S^1$-bundle 
$M$
and let 
$Y$
be a non-empty collection of boundary components of
$M$.
Let 
$S$
be a normal section of 
$M$
and let
$A_1,\ldots,A_n$
be the normal vertical annuli in 
$M$,
that intersect 
$S$
in a disjoint collection of arcs.
We also have:
$S\cap Y=\mathcal{P}$
and
$(A_1\cup\ldots\cup A_n)\cap Y=\mathcal{F}\cap Y$.

The simplified triangulation of the 
$S^1$-bundle
$M$
can be defined as follows. 
Each torus
in
$\partial M -Y$
contains precisely two circles,
one from
$\partial S$
and the other from
$\partial (A_1\cup\ldots\cup A_n)$,
that intersect in a single
point. We can triangulate this torus by 2 triangles.
Any component of
$Y$
inherits a simplicial structure from the triangulation
$T_Y$
and from the families of normal simple closed curves
$\mathcal{P}$
and
$\mathcal{F}$.
We obtain a subdivision of
$T_Y$
by coning each of these PL discs from one 
of the vertices in their boundaries. 
What we defined so far specifies the triangulation of
the boundaries of both our vertical annuli and of the section
$S$.
On each annulus there is precisely one spanning arc coming
from its intersection with the surface 
$S$.
The endpoints of this arc are vertices in the already
defined triangulation of the boundary of the annulus. So 
we can take this arc to be an edge of the simplified 
triangulation. We can then triangulate the vertical annulus
by coning from one of endpoints of the arc.
The arcs from
$S\cap  (A_1\cup\ldots\cup A_n)$
decompose the section 
$S$
into a single disc. The boundary of this disc is already triangulated.
So we can triangulate 
$S$
by coning this disc from some vertex in its boundary. 
Finally, the complement of the polyhedron
$S\cup  A_1\cup\ldots\cup A_n$
in the 
$S^1$-bundle 
$M$,
is a 3-ball. Since the triangulation of its boundary is
already defined, we can take the simplified triangulation
of
$M$
to be the cone in this 3-ball.

Using the simplified triangulation of the manifold
$M$
we can now establish the following key step in the proof
of our main theorem.\\

\begin{thm}
\label{thm:bundles}
Let 
$M\rightarrow B$
be an
$S^1$-bundle over a possibly non-orientable bounded surface
$B$
with negative Euler characteristic. Let 
$Y$
be a non-empty collection of boundary components of the manifold
$M$.
Let 
$P$
and 
$Q$
be two triangulations of
$M$
containing 
$p$
and
$q$
tetrahedra respectively. 
Assume also that the triangulation 
of 
$Y$
induced by
$P$
is isotopic (in
$Y$)
to the triangulation obtained by restricting the triangulation
$Q$
to the components of
$Y$.
Then there exists a sequence of Pachner
moves of length less than 
$e^2(400p^2)+e^2(400q^2)$
that transforms the triangulation
$P$
into a triangulation which is isomorphic to the triangulation
$Q$.
The homeomorphism of
$M$
that realizes this simplicial isomorphism maps 
fibres onto fibres, does not permute
the components of
$\partial M$
and, when restricted to
$Y$,
it is isotopic to the identity on
$Y$.
\end{thm}

\begin{proof}
Instead of transforming the triangulation
$P$
to 
$Q$,
we shall change each one to a simplified triangulation.
It follows directly from lemma~\ref{lem:section} and 
lemma~\ref{lem:twist} that the two simplified triangulations
obtained in this way, can be related by a homeomorphism
that has the required properties. So what we need to find is a bounded
sequence of Pachner moves that will change the triangulation 
$P$
(and
$Q$)
to a simplified triangulation.

We start by fixing a ``minimal'' family 
$\mathcal{P}$
of simple closed curves in the components of
$Y$
that bounds a normal section in
$M$.
Theorem~\ref{thm:section} implies that there
exists a section
$S$
with the prescribed boundary which is normal with respect
to 
$P$,
and is also fundamental.
We now choose a family
$\mathcal{F}$
of fibres in the boundary of 
$M$, 
that consist of the smallest possible number of normal arcs.
Using proposition~\ref{prop:annuli} we have already constructed 
a family of vertical normal annuli
$A_1,\ldots,A_n$
with bounded complexity, whose boundary equals
$\mathcal{F}$.
We also know that we can assume that the intersection between the
annuli and the section consists of spanning arcs only.

We can now apply lemma~\ref{lem:norm} to the triangulation
$P$
and the normal section
$S$.
By proposition~\ref{prop:hass} we have that the number of
normal discs in
$S$
is bounded above by
$5p(7p2^{7p})$.
(The factor
$5p$
comes in since 
$S$
contains at most 5 distinct disc types in any tetrahedron of
$P$.)
So we can subdivide 
$P$
by making less than
$200\cdot35p^32^{7p}$
Pachner moves. The triangulation we get contains the section
$S$
in its 2-skeleton and consists of 
$20(35p^22^{7p}+p)<20\cdot 2^{11p}$
tetrahedra. 
The complementary regions of normal discs of 
$S$
in the tetrahedra of
$P$
are, in this subdivision, triangulated as cones on their boundaries.
So we can assume that the surface
$A_1\cup\ldots\cup A_n$
is still normal in the subdivision. Furthermore, since
each normal disc in
$P$
intersects any 3-simplex of the subdivision in at most one 
disc, we can bound the normal complexity of our vertical
annuli in this new triangulation.

From the discussion in the previous subsection we know that
the number of copies of any disc type in the surface
$A_1\cup\ldots\cup A_n$,
that is contained in a 3-simplex of
$P$,
is bounded above by
$(2\cdot 2^{11p})^{12p}$.
This implies that each normal coordinate of the union of vertical annuli,
in the subdivided triangulation, is certainly bounded above by
$5(2\cdot 2^{11p})^{12p}$.
So the surface
$A_1\cup\ldots\cup A_n$
consists of less than
$20\cdot2^{11p}5(2\cdot 2^{11p})^{12p}$
normal discs
in the subdivision. Applying lemma~\ref{lem:norm} again we obtain
a further subdivision of
$P$
that contains the polyhedron
$S\cup A_1\cup\ldots\cup A_n$
in its 2-skeleton.
This subdivision can be constructed by 
$200(20\cdot 2^{11p})(20\cdot2^{11p}5(2\cdot 2^{11p})^{12p})$
Pachner moves.
The number of 3-simplices in the subdivision is, according to
lemma~\ref{lem:norm}, bounded above by
$20(20\cdot 2^{11p}+20\cdot2^{11p}5(2\cdot 2^{11p})^{12p})<2^{150p^2}$.

We can now apply theorem~\ref{thm:3ball} to the complement of the 
polyhedron
$S\cup A_1\cup\ldots\cup A_n$
in the 3-manifold
$M$.
It gives us that after 
$2^{2^{380p^2}}$
Pachner moves, we can assume that the triangulation of the complement
is a cone on the boundary. The number of tetrahedra in this triangulation
is bounded above by
$2^{150p^2}$.

Our next task is to modify the triangulation of 
both the vertical annuli and the
section
$S$.
The arcs
$S\cap (A_1\cup\ldots\cup A_n)$
decompose all vertical annuli into discs. Since these arcs are 
already contained in the 1-skeleton of the subdivision, we can
apply lemma 4.2 from~\cite{mijatov} to simplify the triangulations
of these discs. Since the discs have a cone on both sides, each
2-dimensional Pachner moves can be implemented by the 3-dimensional
ones in the usual way. This procedure requires linearly many
(in the number of tetrahedra)
Pachner moves. We can thus assume that our triangulation has 
this property already. In a similar way we get the required
triangulation of the section 
$S$.
We also need to amalgamate all 1-simplices contained in the
arcs
$S\cap (A_1\cup\ldots\cup A_n)$.
By lemma 4.1 in~\cite{mijatov} 
we can crush one 1-simplex in an arc by less than a hundred 
Pachner moves. Since there are again linearly many (in the number of
tetrahedra in the subdivision) 1-simplices, we can assume that in
our subdivision the arcs from
$S\cap (A_1\cup\ldots\cup A_n)$
are 1-simplices.

The triangulation we have so far is almost isomorphic to
the simplified triangulation. They differ only on the 
boundary components of 
$M$
that are not contained in 
$Y$.
Let 
$A$
be one such torus. Then a boundary component of 
$S$
and
a boundary circle of some vertical annulus decompose 
$A$
into a disc. We can change the triangulation of this disc into
a cone on its boundary using 
lemma 4.2 in~\cite{mijatov}. We can implement the
2-dimensional Pachner moves by the 3-dimensional ones. Here we have 
to use the 3-dimensional moves that alter the triangulation
of the boundary. All this can be done by linearly many moves.
Also, we have to change the triangulation of the boundary curves
of 
$S$
and the vertical annulus in
$A$.
This too can be achieved by linearly many Pachner moves. This 
finally brings
us to the simplified triangulation.

Putting everything together we can conclude that the whole 
process that is described above takes not more than
$e^2(400p^2)$
Pachner moves. This proves the theorem.
\end{proof}

%% file: exceptional.tex
We are now going to consider the triangulations of the 
following 3-manifolds:
$S^1\times S^1\times I$,
the unique orientable
$S^1$-bundle over the Moebius band,
the doubles of these two manifolds, and
a Seifert fibred space
over a projective plane with two singular
fibres with indices
$\frac{1}{2}$
and
$-\frac{1}{2}$.
The reason why our usual 
techniques fail in this context is because these manifolds contain
horizontal surfaces of zero Euler characteristic. 

Let's first deal with the two product manifolds
$S^1\times S^1\times I$
and
$S^1\times S^1\times S^1$.
Our original strategy was to find essential surfaces that
are vertical in a given fibration of the manifold. In 
these cases 
the spaces
fibre in infinitely many ways, but all of these fibrations
are 
equivalent under a homeomorphism. 
Loosely speaking what we
do now is precisely the opposite from before. We first choose an essential
surface and then decide which fibration we are going to use. 

Assume that at least one component of 
the boundary of 
$S^1\times S^1\times I$
has a prescribed simplicial structure. 
What this means is that we are not allowed to use homeomorphisms
of the manifold which are not isotopic to the identity on the
boundary.
We fix an isotopy class 
which contains a simple closed curve in the 
boundary component of the product
$S^1\times S^1\times I$,
that consists of the smallest number of normal arcs among all
homotopically non-trivial curves.
If both components have a prescribed triangulation, we first choose a
component and then fix the curve in it. 
Let
$S$
be a minimal weight normal annulus among all annuli with one
boundary circle in the chosen isotopy class and the other
circle in the opposite component of
$\partial (S^1\times S^1\times I)$.
If
$S=F+G$,
then by theorem~\ref{thm:sum} we can conclude that 
the surfaces
$F$
and
$G$
are both incompressible and 
$\partial$-incompressible. Since their Euler characteristics
have to vanish (there are no trivial patches by lemma~\ref{lem:patches})
and since
$S^1\times S^1\times I$
contains no embedded Moebius bands, at least one of them, say 
$F$,
has to be 
an annulus. 
It follows immediately that 
$F$ 
has a possibly ``shorter'' boundary than
$S$.
But its weight is certainly smaller than that of
$S$.
This contradiction implies that 
$S$
has to be fundamental. 

Now we can use the vertical annulus 
$S$
to construct a vertical solid torus
in the interior of the manifold, without changing the triangulation
of the boundary. 
This is achieved by taking a thin regular neighbourhood of our
annulus that stops just before the boundary of the manifold and then
using the same technique as in the proof of
proposition~\ref{prop:prod}. After 
removing this solid torus from the manifold, the 
Euler characteristic of the base surface becomes negative
and we are back in the ``generic'' case.

If the triangulation of the boundary is not fixed, we can still 
find a fundamental vertical annulus as described above. 
In the product manifold
$S^1\times S^1\times I$
we can always find a homeomorphism mapping one choice to the other.
This proves theorem~\ref{thm:main}
in this exceptional case.

It is well known that any two incompressible tori in
$S^1\times S^1\times S^1$
are related by a homeomorphism of the ambient manifold. 
So we are free
to pick the one which has smallest weight. Such a torus is
fundamental by theorem~\ref{thm:sum}, because the manifold
contains no embedded Klein bottles. We then look for 
an essential annulus in the complement of the regular neighbourhood
of the torus. This can be done in the same way as above. Again
we apply the techniques of proposition~\ref{prop:prod} to obtain
a vertical solid torus neighbourhood of a regular fibre. After
we remove it we can use the same strategy as in the ``generic'' case
to further simplify the triangulation. 

In case of the 
$S^1$-bundles over a Moebius band 
we are going to use the fact that this manifold is homeomorphic
to a Seifert fibred space over a disc with two exceptional fibres
of index   
$\frac{1}{2}$.
This situation was dealt with in section~\ref{sec:detecting}
when we were isolating singular fibres of such Seifert 
fibred spaces (c.f. the exceptional case in the construction
of family
$\mathcal{B}$).
In case our 
$S^1$-bundle has a prescribed triangulation of its boundary,
a combination of the strategy for
$S^1\times S^1\times I$
and the techniques mentioned above prove the main theorem.

The second homology group
$H_2(M;\ZZ)$,
where 
$M$
is our 
$S^1$-bundle over a Klein bottle,
is isomorphic to 
$\ZZ$.
The vertical torus above the non-separating orientation
preserving simple closed curve in the Klein bottle
is a generator of
$H_2(M;\ZZ)$.
The results from~\cite{wang} imply that
the least weight normal representative
in this homology class has to be a fundamental torus.
Furthermore this fundamental torus is isotopic to the vertical
one we started with. Its complement 
is therefore homeomorphic to the product
$S^1\times S^1\times I$.
So 
$M$
supports a torus bundle structure with holonomy equal to
$-\mathrm{id}$.
Notice that any homeomorphism of the fibre
$S^1\times S^1$
extends to the whole of
$M$.
We can thus choose a spanning annulus of minimal weight in the 
complement 
$S^1\times S^1\times I$
which is vertical up to a homeomorphism
and which is also fundamental. Removing its regular
neighbourhood yields an 
$S^1$-bundle over a bounded surface and our usual techniques 
apply.

Let 
$M$
be
a Seifert fibred space
over a projective plane with two singular
fibres with indices
$\frac{1}{2}$
and
$-\frac{1}{2}$.
There is only one essential vertical torus 
$A$
in
$M$
up to homeomorphism.
It is lying above the simple closed curve in
$\RR P^2$
that cuts off a disc containing both singular points.
Since 
$H_2(M,\ZZ)$
is trivial, any horizontal torus 
$A'$
has to be separating.
Both components 
$X$
and
$Y$
of the complement are homeomorphic to
an 
$I$-bundle over a Klein bottle. Any incompressible 
$\partial$-incompressible
annulus in 
$X$
and
$Y$
is vertical in some Seifert fibration of the pieces.
The torus 
$A$
gives us at least two such annuli, one in
$X$
and one in 
$Y$.
Their respective boundary slopes are identified
by the gluing map from
$\partial X$
to
$\partial Y$.
This
means that these fibrations form a new Seifert 
fibration of
$M$
in which 
$A'$
is vertical. 
But by the classification of Seifert fiberings it follows
that 
$M$
has a unique Seifert fibration up to homeomorphism.  
So there is a homeomorphism of
$M$
that maps 
$A'$
onto
$A$.
In other words we have proved that there is only one
incompressible torus
in
$M$
up to homeomorphism. It also follows immediately that 
$M$
contains precisely one injective Klein bottle up 
to homeomorphism (the boundary of its regular neighbourhood
is the unique incompressible torus).
We now pick a normal surface which has the smallest weight
out of all essential tori and injective Klein bottles in
$M$.
By theorem 4.1 from~\cite{bart} we can conclude that such a
surface is fundamental. It is clear now how to construct
an essential torus in
$M$. 
We can then use the same method as in 
$S^1\times S^1\times S^1$
to finish the proof of theorem~\ref{thm:main} for the 
exceptional manifolds.

%% file: conc.tex
We have now developed all the necessary tools to prove theorem~\ref{thm:main}.
Section~\ref{sec:exceptional} deals with the case when the space 
$M$
is homeomorphic to one of the exceptional manifolds. The 
strategy in the generic case is as follows. We start by applying 
theorem~\ref{thm:T2}
to both triangulations 
$P$
and
$Q$
of the Seifert fibred space 
$M$.
If we are in the situation where there are no singular fibres,
we apply proposition~\ref{prop:prod} instead. This produces
a subdivision 
$P_2$ 
of the triangulation
$P$
containing less than
$e^3(500p^2)$
tetrahedra. The process requires less than
$e^4(500p^2)$
Pachner moves. The subdivision
$Q_2$,
obtained in the similar way from triangulation
$Q$,
induces a triangulation in
the components of
$\partial M$
that belong to 
$Y$,
which is isotopic to the one induced by
$P_2$.
Also both triangulations
$P_2$
and
$Q_2$
contain subcomplexes that triangulate 
the same 
$S^1$-bundle with non-empty boundary.
This bundle is obtained by removing the nicely
triangulated regular neighbourhoods of 
the singular (or possibly one regular) fibres in
$M$.
It follows from the definition of the subdivisions
$P_2$
and
$Q_2$
that the simplicial structures in the newly acquired boundary 
components coincide. So we can add them to the collection
$Y$.

Theorem~\ref{thm:bundles} can now be applied. It implies that we
can transform one of the subcomplexes, triangulating
the
$S^1$-bundle,
into another by making at most
$e^2(400e^3(500p^2)^2)+ e^2(400e^3(500q^2)^2)$
Pachner moves.
It is also clear that the homeomorphism of the bundle, described in 
theorem~\ref{thm:bundles},
extends to the whole 3-manifold
$M$.
It obviously has all the required properties that are stated in
theorem~\ref{thm:main}.
Putting all the above numbers together concludes the proof.